\magnification \magstep1
\input amstex
\documentstyle{amsppt}
\rightheadtext{Normal forms for hypersurfaces}
\topmatter
\leftheadtext{Peter Ebenfelt}
\date{\number\year-\number\month-\number\day}\enddate


\define \im{\text{\rm Im }}
\define \re{\text{\rm Re }}

\define \bR{\Bbb R}

\define \bC{\Bbb C}



\def\dim {\text {\rm dim}}

\def\Qb{\overline{Q}}
\def\fb{\overline{f}}
\def\gb{\overline{g}}
\def\fob{\overline{f^1}}
\def\ftb{\overline{f^2}}
\def\Aut{\text{\rm Aut}}
\def\pb11{\overline{p_{11}}}
\def\pb02{\overline{p_{02}}}
\def\sign{\text{\rm sgn }}

\title Normal forms and biholomorphic equivalence of real hypersurfaces
in $\bC^3$
\endtitle
\author Peter Ebenfelt\footnote{Supported in part by a grant from the
Swedish Natural Science Research
Council.\hfill\hfill}
\endauthor
\affil Department of Mathematics\\Royal Institute of Technology\\100
44 Stockholm, Sweden\endaffil
\address Department of Mathematics, Royal Institute of Technology, 100
44 Stockholm, Sweden\endaddress
\email ebenfelt\@math.kth.se\endemail

\toc
\widestnumber\head{3.3.3}
\head 1. Introduction\endhead
\head 2. Finite nondegeneracy and holomorphic nondegeneracy\endhead
\head 3. Normal forms and biholomorphic equivalence -- the main results\endhead
\head 4. Some applications and examples\endhead
\subhead 4.1. Computation of two fourth order invariants\endsubhead
\subhead 4.2. Everywhere Levi degenerate hypersurfaces\endsubhead
\head 5. Regular coordinates and preliminaries\endhead
\head 6. Proof of Theorem A (i)\endhead
\head 7. Proof of Theorem A (ii); the beginning\endhead
\subhead 7.1. The setup\endsubhead
\subhead 7.2. The case $\Delta_{bc}=\Delta_{ac}=0$\endsubhead
\subhead 7.3. The case $\Delta_{ab}=\Delta_{ac}=0$\endsubhead
\subhead 7.4. The case $\Delta_{ab}=\Delta_{bc}=0$\endsubhead
\head 8. Proof of Theorem A (ii); the conclusion\endhead
\subhead 8.1. Two simple lemmas\endsubhead
\subhead 8.2. The case $\Delta_{ac}=0$ and $\Delta_{ab}\Delta_{bc}\neq0$\endsubhead
\subhead 8.3. The case $\Delta_{ab}=0$ and $\Delta_{ac}\Delta_{bc}\neq0$\endsubhead
\subhead 8.4. The case $\Delta_{bc}=0$ and $\Delta_{ab}\Delta_{ac}\neq0$\endsubhead
\subhead 8.5. Last case; $\Delta_{ab}\Delta_{ac}\Delta_{bc}\neq0$\endsubhead
\head 9. Proof of Theorem B\endhead
\subhead 9.1. The setup\endsubhead
\subhead 9.2. Conclusion of the proof of Lemma 9.1.18; the cases
(A.i.1) and (A.i.3)\endsubhead
\subhead 9.3. Conclusion of the proof of Lemma 9.1.18; the case
(A.i.2)\endsubhead
\endtoc

\endtopmatter
\document

\heading 1. Introduction\endheading

As Poincar\'e
[P] perhaps was the first to realize, 
a real hypersurface $M$ in $\bC^N$,  when $N\geq2$, has  non-trivial
local invariants under biholomorphic transformations of $\bC^N$ near
a distinguished point $p_0\in M$. To understand the biholomorphic
equivalence class of $(M,p_0)$, i.e. the class consisting of those
real hypersurfaces $M'\subset\bC^N$ with distinguished points $p_0'\in
M'$ that are locally equivalent to $(M,p_0)$ under some germ of a
biholomorphic transformation, is a fundamental problem.
It is well 
understood in the case where $M$ is real analytic and Levi nondegenerate at
$p_0$. The solutions are due to Cartan [C1--C2] in the case $N=2$, and
to Chern and Moser [CM]
in the case of general $N$. Important contributions were also
made by Tanaka [T1--T2]; see also [BS]. 

In contrast to the situation where $M$ is Levi nondegenerate at $p_0$,
very little is known at Levi
degenerate points. We should mention though that Webster, in [W], studied
the equivalence under holomorphic contact
transformations (a larger class of transformations than that of
biholomorphic ones) of real hypersurfaces at Levi degenerate points; under
such transformations, e.g.,  all Levi nondegenerate real hypersurfaces are
equivalent. In the present
paper, we shall consider real-analytic
hypersurfaces in
$\bC^3$ at Levi degenerate points of a certain kind. Together with the results of Chern-Moser, 
the results obtained here yield a fairly complete understanding of the
biholomorphic equivalence problem at generic points on real-analytic
hypersurfaces in $\bC^3$.

In [CM], two solutions of the biholomorphic equivalence problem for
Levi nondegenerate hypersurfaces are
given: an intrinsic solution in terms of classical objects of
differential geometry and an extrinsic solution in terms of a normal
form for $M$ at $p_0$. The extrinsic approach, which is also the
approach taken 
in this paper, amounts to explicitly defining a class of real
hypersurfaces $\Cal N$ at $0\in\bC^N$ such that $M$ can be
transformed to one of the surfaces in $\Cal N$ by a biholomorphic
transformation $Z'=H(Z)$ near $p_0$. Moreover, this transformation
should be
unique, given that a finite dimensional choice of normalization
conditions (which in the Levi nondegenerate case can be identified
with a choice of element in a classical group of rational
transformations) for $H$ are made. Having defined such a normal form,
one obtains the following answer to
the biholomorphic equivalence problem: 

{\it Two germs $(M,p_0)$ and
$(M',p'_0)$ of hypersurfaces are biholomorphically equivalent
(i.e.\ there is a germ at $p_0$ 
of a biholomorphic transformation $Z'=H(Z)$ such that $H(p_0)=p_0'$
and $H(M)\subset M'$) if and only if for two
(possibly different) choices of normalizations they can be brought to
the same normal form.}

The results that we obtain in this paper can be described as
follows. We consider the 
class of real-analytic hypersurfaces $M$ in $\bC^3$ that are
2-nondegenerate (see \S2 for the definition) at a distinguished point
$p_0$. We first obtain a partial normal form for $M$ at $p_0$ (Theorem
A in \S3) prescribing the lowest order terms (i.e.\ the data associated with
2-nondegeneracy) in the defining equation of $M$. This divides the
class of 2-nondegenerate hypersurfaces into 8 different types, called (A.i.1--3),
(A.ii.1--5) in Theorem A. Moreover, the types (A.ii.1) and
(A.ii.3--5) come
with one or two non-trivial invariants. 
We then obtain
a complete, but formal, normal form, i.e.\ we obtain a class of formal
hypersurfaces 
such that $M$ can be transformed into one of these by a formal
invertible transformation which is unique with a finite dimensional
choice of normalization, for three of these types (A.i.1--3); see Theorem B in \S3. 
This
implies that any two germs $(M,p_0)$ and $(M',p_0')$ of types
(A.i.1--3) are formally equivalent if and only if for two choices of
normalizations the two can be brought to the same formal normal form.
Now, in view
of a theorem from [BER3] (see also \S2), any formal invertible transformation taking
a real-analytic, finitely nondegenerate hypersurface $M$ into another
such hypersurface $M'$ is biholomorphic and, hence, formal
equivalence of finitely nondegenerate hypersurfaces is the same as
biholomorphic equivalence.  

As we mentioned above, these
results together with the results 
from [CM] give a fairly complete understanding of the biholomorphic
equivalence problem for real-analytic hypersurfaces in $\bC^3$ at
generic points, because a real-analytic hypersurface $M$ in $\bC^3$ is
at a generic point $p\in M$ either (a) locally of the form $\tilde
M\times\bC$, $\tilde M\subset\bC^2$, in a
neighborhood of $p$, (b) 2-nondegenerate with one non-zero eigenvalue
of the Levi form, or (c) Levi nondegenerate. (This will be further explained
in \S2). If $M$ is foliated by
complex curves as described in (a), then one may reduce the biholomorphic equivalence
problem to one for hypersurfaces in $\bC^2$ where the problem is well
understood (again, at generic points). If $M$ is 2-nondegenerate at
$p\in M$ with
one non-zero eigenvalue of the Levi form at that point, then $M$ is of
one of the types (A.i.1--3); in fact, we show (Theorem 4.2.8) that if
$M$ is 2-nondegenerate at {\it generic points} (which implies that $M$ is
everywhere Levi degenerate), then $M$ is of type (A.i.2) at every
2-nondegenerate point. Finally, if $M$ is Levi
nondegenerate at $p\in M$, then the Chern--Moser theory applies.

In \S4, we discuss a few applications of the results in \S3. We
compute two invariants of fourth order for a real hypersurface of type
(A.i.2) and consider holomorphically nondegenerate, but everywhere Levi
degenerate, real-analytic hypersurfaces. We also discuss a
few examples in this section, such as e.g. 
Freeman's example [F] of an everywhere Levi degenerate hypersurface
which is not foliated by complex curves.


\heading 2. Finite nondegeneracy and holomorphic
nondegeneracy\endheading  

Let $M$ be a real-analytic hypersurface in $\bC^N$ and let $p_0\in
M$. Following Stanton [S], we say
that $M$ is {\it holomorphically nondegenerate} at $p_0$ if there is
no germ at $p_0$ of a holomorphic vector field
$$
X=\sum_{k=1}^Na_k(Z)\frac{\partial}{\partial Z_k},\tag 2.1
$$
where the $a_k(Z)$ are germs at $p_0$ of holomorphic functions, such
that $X$ is
tangent to $M$ near $p_0$. If $\rho(Z,\bar Z)=0$ is a local defining
equation for $M$ near $p_0$, then we say that $M$ is {\it finitely
nondegenerate} at $p_0$ if there is a positive integer $k$ such that
$$
\text{\rm span}\,\{L^\alpha\rho_Z(p_0,\bar p_0)\:|\alpha|\leq
k\}=\bC^N;\tag 2.2
$$
here, $L_1,\ldots,L_{N-1}$ is a basis for the CR vector fields of $M$
near $p_0$, 
$$
L^\alpha=L_1^{\alpha_1}\ldots L_{N-1}^{\alpha_{N-1}}, 
$$
and $\rho_Z$ denotes the holomorphic gradient $\partial\rho/\partial
Z$. This notion is independent of the choice of defining function, the
choice of basis for the CR vector fields, and the coordinates used
([BHR]). Moreover, the smallest $k$ for which \thetag{2.2} holds is a
biholomorphic (and formal) invariant and we say, more precisely, that $M$ is
$k$-nondegenerate at $p_0$ if $k$ is the smallest integer for which
\thetag{2.2} holds. As is easy to verify, $M$ is Levi nondegenerate at
$p_0$ if and only if $M$ is 1-nondegenerate at $p_0$. 

Some of the main facts
(whose proofs can be found in e.g. [BER1], where generalizations to
higher codimensional CR manifolds also can be found)
that we shall need about these notions can be summarized in the
following proposition. 

\proclaim{Proposition 2.3} Let $M\subset\bC^N$ be a connected
real-analytic hypersurface. The following are
equivalent. 
\roster 
\item"(i)" There exists $p_1\in M$ such that $M$ is holomorphically
nondegenerate at $p_1$. 
\item"(ii)" $M$ is holomorphically nondegenerate at every point $p\in M$.
\item"(iii)" There exists $p_1\in M$ such that $M$ is finitely
nondegenerate at $p_1$. 
\item"(iv)" There exists a proper real-analytic subset $V$ of $M$ and
an integer $\ell=\ell(M)$, with $1\leq \ell\leq N-1$, such that $M$ is 
$\ell$-nondegenerate at every $p\in M\setminus V$. 
\endroster
\endproclaim

We say that a connected real-analytic hypersurface is
holomorphically nondegenerate if it is so at one point (and hence at
all points). If $M$ is holomorphically nondegenerate, then the number
$\ell(M)$ provided by Proposition 2.3 (iv) is called the {\it Levi
number} of $M$. 

Before we turn to a discussion of real-analytic hypersurfaces in
$\bC^3$, we wish to state two general results concerning transformations between
finitely nondegenerate hypersurfaces in $\bC^N$. The first result
(from [BER3]) will not be used in this paper, but it explains why
finitely nondegenerate hypersurfaces seem to the right class of
hypersurfaces for which one can hope to obtain a normal form. 
Recall that we would like the transformation of a germ $(M,p_0)$ to
normal form to be 
essentially unique (at least up to a finite dimensional
normalization), and a transformation up to normal form can only be
unique up to transformations preserving $(M,p_0)$ (the stability
group). It is not difficult to see that the stability group of a
holomorphically degenerate hypersurface at a generic point is infinite
dimensional (in fact, this is true at every point; see [BER2]). In
contrast, we have the following result from [BER3] (cf. also [BER2])
for finitely nondegenerate hypersurfaces. 
We denote by $\Aut(M,p_0)$ the stability group of $M$ at $p_0$,
i.e. the group of germs at $p_0$ of biholomorphic transformations
preserving the germ $(M,p_0)$, equipped with its natural inductive
limit topology.

\proclaim{Theorem 2.4 ([BER3])} Let $M\subset\bC^N$ be a real-analytic
hypersurface that is finitely nondegenerate at $p_0\in M$. Then the
stability group $\Aut(M,p_0)$ is a finite dimensional Lie
group.\endproclaim 

\flushpar
{\it Remark.} The arguments in [BER2--3] also give a
bound on the dimension of $\Aut(M,p_0)$ in terms of the complex dimension $N$
of the ambient space and $k$, the order of finite nondegeneracy of $M$
at $p_0$. In the special case of a 2-nondegenerate hypersurface in
$\bC^3$, this bound is $\dim_\bR\Aut(M,p_0)\leq 102$. A consequence of
Theorem B of the present paper is an improved bound on this
dimension in case $M$ is of either of the types (A.i.1--3) as
described by Theorem A.\medskip  

Another result, which will be important in this paper, is the following
result from [BER3] 
that, as we explained in the introduction, reduces the problem of
biholomorphic equivalence of two germs $(M,p_0)$ and $(M',p_0')$ to
that of formal equivalence. We say that an $N$-tuple of formal power
series in $Z$, $H(Z)=(H_1(Z),...,H_N(Z))$, is a formal equivalence
between $(M,0)$ and $(M',0)$ in $\bC^N$ (with the obvious generalization to
general points $p_0\in M$ and $p_0'\in M'$) if $H(Z)$ has no constant
term, $H(Z)$ is invertible (i.e.\ the linear part of $H(Z)$ is
invertible), and if there is a formal power series $a(Z,\bar Z)$ such
that 
$$
\rho'\left(H(Z),\overline{H(Z)}\right)=a(Z,\bar Z)\rho(Z,\bar
Z),\tag2.5
$$
where $\rho(Z,\bar Z)=0$ is a defining equation for $M$ at $0$ and
$\rho'(Z,\bar Z)=0$ is a defining equation of $M'$ at $0$.

\proclaim{Theorem 2.6 ([BER3])} Let $M$ and $M'$ be real analytic
hypersurfaces in $\bC^N$ that are finitely nondegenerate at $p_0\in
M$ and $p_0'\in M'$, respectively. If $H(Z)$ is a formal
equivalence between $(M,p_0)$ and $(M',p_0')$, then $H(Z)$ is
convergent i.e.\ there is a biholomorphic equivalence between
$(M,p_0)$ and $(M',p_0')$ whose power series coincides with
$H(Z)$.\endproclaim 

Let us conclude this section with a discussion of
hypersurfaces in $\bC^2$ and $\bC^3$. In $\bC^2$, a real-analytic
hypersurface is either holomorphically 
degenerate or Levi nondegenerate at all points outside a proper
real-analytic subset. If $M\subset \bC^2$ is holomorphically
degenerate, then it is in fact locally biholomorphic to a real
hyperplane (this is not difficult to see). Thus, in $\bC^2$ every
real-analytic hypersurface is 
either biholomorphic to a real hyperplane or Levi-nondegenerate at
generic points (and thus covered by the
Cartan--Chern--Moser theory at such points).

However, already in $\bC^3$ the situation is a little more
complicated. If $M\subset \bC^3$ is holomorphically degenerate, then
at every point $p_0\in M$ there is a holomorphic vector field $X$ of the
form 
\thetag{2.1} near $p_0$ that is tangent to $M$. Outside a proper
real-analytic subset $V$ of $M$ near $p_0$, the vector field $X$ is
nonsingular (i.e.\ at least one of the coefficients is non-zero). At
such a point $p_1\in M\setminus V$, it follows from the Frobenius
theorem that $M$ is foliated by complex curves and it is possible to choose
coordinates near $p_1$ such that $M$ is given locally as $\tilde M\times \bC$,
for some $\tilde M\subset \bC^2$. Thus, at generic points on a
holomorphically degenerate hypersurface in $\bC^3$ the biholomorphic
equivalence problem may be reduced to the problem in $\bC^2$. 
If $M\subset\bC^3$ is holomorphically nondegenerate, then it is at
generic points either Levi nondegenerate, in which case it is covered
by Chern--Moser theory, or it is $2$-nondegenerate. Hence, there is a
gap between those hypersurfaces that can be studied 
by Chern--Moser theory at generic points and those for which the biholomorphic
equivalence problem at generic points can be reduced to a lower dimensional
problem. This observation was the original motivation for the present work.

As will be further discussed in the following sections, a
real-analytic hypersurface $M$ in $\bC^3$ can {\it a priori} be
$2$-nondegenerate at a point $p_0\in M$ in
two different ways: the Levi form of $M$ at $p_0$ has (i) one
zero and one non-zero eigenvalue or (ii) both eigenvalues of the Levi form
are zero. A moments reflection will convince the reader that outside a
proper real-analytic subset of $M$, the Levi form has at least one
non-zero eigenvalue because if both eigenvalues of the Levi form were
zero on an open subset of $M$ then $M$ would be locally biholomorphic to a
real hyperplane and thus, in particular, holomorphically
degenerate. A closer analysis shows that the situation (i) occurs in
essentially three different ways (A.i.1--3), and the situation (ii) in
5 different ways (A.ii.1--5), as we shall see in Theorem A. On the
other hand, as we shall prove in Theorem 4.2.8, if $M$ is everywhere
Levi degenerate, then at a 2-nondegenerate point it must be of type (A.i.2).

\heading 3. Normal forms and biholomorphic equivalence -- the main results\endheading

Let $M$ be a real-analytic hypersurface in $\bC^3$ and let $p_0\in M$ be a
point at which $M$ is $2$-nondegenerate. We may choose coordinates
$(z,w)=(z_1,z_2,w)\in \bC^3$, vanishing at $p_0$, such that $M$ is
described locally near $p_0=0$ by the equation
$$
\im w=\phi(z,\bar z,\re w),
$$
where $\phi(z,\bar z,s)$ is a real-analytic function with
$\phi(0,0,0)=0$ and $d\phi(0,0,0)\neq 0$. 
Our first result describes a partial normal form for $(M, p_0)$. 

\proclaim{Theorem A} Let $M$ be a real-analytic hypersurface in
$\bC^3$ and assume that $M$ is $2$-nondegenerate at $p_0\in M$. Then
$(M, p_0)$ is biholomorphically 
equivalent to $(M',0)$, where $M'$ is a real-analytic hypersurface of
one the following model forms. 
\medskip
{\rm (i)} If the Levi form of $M$ at $p_0$ has precisely one non-zero
eigenvalue, then $M'$ is one of the following:
$$
\im w=|z_1|^2+|z_2|^2(z_2+\bar z_2)+\gamma(z_1^2\bar z_2+\bar
z_1^2z_2)+O(|z|^4+|\re w||z|^2),\tag A.i.1
$$
where $\gamma=0,1$;
$$
\im w=|z_1|^2+(z_1^2\bar z_2+\bar
z_1^2z_2)+O(|z|^4+|\re w||z|^2);\tag A.i.2
$$
$$
\im w=|z_1|^2+|z_2|^2(z_1+\bar z_1)+O(|z|^4+|\re w||z|^2).\tag A.i.3
$$
\medskip
{\rm (ii)} If the Levi form of $M$ at $p_0$ is $0$, i.e. both
eigenvalues of the Levi form are zero, then $M'$ is one of the following:
$$
\im w=|z_1|^2(z_2+\bar z_2)+r(z_1^2\bar z_2+\bar z_1^2z_2)+
O(|z|^4+|\re w||z|^2),\tag A.ii.1
$$
where $r>0$;
$$
\im w=|z_1|^2(z_2+\bar z_2)+(z_1^2\bar z_2+\bar
z_1^2z_2)+i|z_1|^2(z_1-\bar z_1)+
O(|z|^4+|\re w||z|^2);\tag A.ii.2
$$
$$
\im w=|z_1|^2(z_2+\bar z_2)+(z_1\bar z_2^2+\bar
z_1z_2^2)+|z_2|^2(\lambda z_2+\bar \lambda \bar z_2)+
O(|z|^4+|\re w||z|^2),\tag A.ii.3
$$
where $\lambda\in\bC$, $\lambda\neq0$;
$$
\aligned
\im w= &|z_1|^2(z_1+\bar z_1)+|z_2|^2(z_2+\bar z_2)+(\mu z_1^2\bar z_2+
\bar \mu\bar
z_1^2z_2)+\\&(\nu z_1\bar z_2^2+\bar \nu\bar
z_1z_2^2)+
O(|z|^4+|\re w||z|^2),\endaligned\tag A.ii.4
$$
where $\mu,\nu\in\bC$, $\mu\nu\neq1$.
$$
\aligned
\im w= &|z_1|^2(\eta z_1+\bar\eta\bar z_1)+(z_1^2\bar z_2+
\bar z_1^2z_2)+\\&(z_1\bar z_2^2+\bar
z_1z_2^2)+
O(|z|^4+|\re w||z|^2),\endaligned\tag A.ii.5
$$
where $\eta\in\bC$.
\medskip
\flushpar
Moreover, all of these models can be taken in regular form (see
the next section) and then they are mutually
non-equivalent, provided that we in \thetag{A.ii.4} arrange so that
$|\mu|\geq |\nu|$ and $\arg \mu\geq\arg \nu$, where
$\arg\mu,\arg\nu\in[0,2\pi)$, if $|\mu|=|\nu|$. 
\endproclaim 

\flushpar
{\it Remark.} The models described in Theorem A need not be in regular
form (as described in \S4) in order to be mutually non-equivalent. In
fact, it suffices that 
the remainder in (A.i.1--3) is $O(|z|^4+|\re w||z|^2+|\re w|^2)$, or
$O(4)$ in the weighted coordinate system where $z$, $\bar z$ have
weight one and $\re w$ has weight two, and that the remainder in
(A.ii.1--5) is $O(|z|^4+|\re w||z|+|\re w|^2)$, or $O(4)$ in the
weighted coordinate system where $z$, $\bar z$ have weight one and $\re
w$ has weight three. This follows from the fact that if the equation
for $M$ is in regular form modulo terms of
weighted degree $\nu$, then the entire equation
may be transformed to regular form without changing the terms of
weighted degree less than $\nu$. We refer the reader to
the forthcoming book [BER4] for a proof 
of this.\medskip

The proof of Theorem A will be given in \S5--8.
As our second result, we present a complete formal normal form for $(M, p_0)$ of types
(A.i.1--3) above. In order to describe this result, we need to
introduce some notation. 
We subject a germ $(M,0)$, of either of the types (A.i.1--3), to a formal
invertible transformation
$$
z=\tilde f(z',w')\quad,\quad w=\tilde g(z',w'),\tag3.1
$$
where $\tilde f=(\tilde f^1,\tilde f^2)$, such that the form (A.i.1--3) is preserved. 
We assign the weight one to the variables $z=(z_1,z_2)$, the weight
two to $w$, and say that a polynomial $p_\nu(z,w)$ is weighted
homogeneous of degree $\nu$ if, for all $t>0$,
$$
p_\nu(tz,t^2w)=t^\nu p_\nu(z,w).\tag3.2
$$
We shall write $O(\nu)$ for terms of weighted degree greater than or equal to
$\nu$. Similarly, we speak of weighted homogeneity of degree $\nu$ and
$O(\nu)$ for polynomials and power series in $(z,\bar z,\re w)$, where
$\bar z$ is assigned the weight one and $\re w$ the weight two.
A detailed inspection of the proof of Theorem A (i) (in \S6) yields the following.
We leave the details of this to the reader. 

\proclaim{Proposition 3.3} A transformation \thetag{3.1} preserving
regular coordinates (see the next section) also preserves the
form {\rm (A.i.1--3)} if and only if:
\roster
\item"(1)" In the case {\rm (A.i.1)} with $\gamma=1$, the mapping is of the form
$$
\aligned
\tilde f^1(z,w)& =z_1+Dw-2i\bar Dz_1^2+O(3)\\
\tilde f^2(z,w)& =z_2+O(2)\\
\tilde g(z,w)& =w+2i\bar Dz_1w+O(4),\endaligned\tag3.4
$$
where $D\in\bC$.
\item"(1')" In the case {\rm(A.i.1)} with $\gamma=0$, the mapping is
of the form
$$
\aligned
\tilde f^1(z,w)& =C^{1/2}e^{i\theta}z_1+Dw-2i\bar De^{2i\theta}z_1^2+O(3)\\
\tilde f^2(z,w)& =C^{1/3}z_2+O(2)\\
\tilde g(z,w)& =Cw+2i\bar DC^{1/2}e^{i\theta}z_1w+O(4),\endaligned\tag3.4'
$$
where $C>0$, $\theta\in \bR$, $D\in\bC$.
\item "(2)" In the case {\rm (A.i.2)}, the mapping is of the form
$$
\aligned
\tilde f^1(z,w)& =C^{1/2}e^{i\theta}z_1+Dw-e^{2i\theta}(2i\bar D+C^{1/2}\bar Ae^{i\theta})z_1^2+O(3)\\
\tilde f^2(z,w)& =Az_1+e^{2i\theta}z_2+O(2)\\
\tilde g(z,w)& = Cw+2i\bar DC^{1/2}e^{i\theta}z_1w+O(4)
,\endaligned\tag3.5
$$
where $C>0$, $\theta\in\bR$, $A,D\in \bC$.
\item "(3)" In the case {\rm (A.i.3)}, the mapping is of the form
$$
\aligned
\tilde f^1(z,w)& =C^{1/2}z_1+Dw-2i\bar Dz_1^2+O(3)\\
\tilde f^2(z,w)& =C^{1/4}e^{i\theta}z_2+O(2)\\
\tilde g(z,w)&= Cw+2i\bar D C^{1/2}z_1w+O(4),\endaligned\tag3.6
$$
where $C>0$, $\theta\in\bR$, $D\in\bC$.
\endroster
\endproclaim
 
We shall consider formal mappings \thetag{3.1} of the following form
$$
(\tilde f(z,w),\tilde g(z,w))=(T\circ P)(z,w).\tag 3.7
$$
Here, $P(z,w)$ is a polynomial mapping 
$$
P(z,w)=(A^1_1z_1+Dw+Bz_1^2+q_1(z,w),A^2_1z_1+A^2_2z_2+q_2(z,w),Cw+Ez_1w),\tag3.8
$$
where $A^1_1,A^2_1,A^2_2,B,D,E\in\bC$, $C>0$, are such that $P(z,w)$
satisfies the conditions (depending on the form of $(M,0)$) imposed by
Proposition 3.3, $q_1$, 
$q_2$ are weighted homogeneous polynomials such that $q_1$ is
$O(3)$ and  $q_2$ is $O(2)$, and $T(z,w)$ is a formal mapping
$$
T(z,w)=(z+f(z,w),w+g(z,w)),\tag3.9
$$
where $f=(f^1,f^2)$, and $g$ are formal power series in $(z,w)$ such
that $f^1$ is $O(3)$, $f^2$ is $O(2)$, and $g$ is $O(4)$. We shall
require that the
polynomials $q_1$, $q_2$, and the formal 
series $f^1$, $f^2$, $g$ satisfy additional conditions that will be
different depending on
which of the types (A.i.1--3) the germ $(M,0)$ is. 

If $(M,0)$ is of the form (A.i.1) or (A.i.3), then we shall require
that the polynomials $q_1$, $q_2$ in
\thetag{3.8} are of the form
$$
q_1(z,w)=Rz_1w+\sum_{|\beta|=3} C_\beta z^\beta,\quad
q_2(z,w)=\sum_{|\alpha|=2}D_\alpha z^\alpha,\tag3.10
$$
for $C_\beta,D_\alpha\in\bC$ and $R\in\bR$, and that the formal
series $f^1$, $f^2$ are such that the constant terms in the following
formal series vanish (the index $k$ below ranges over $\{1,2\}$,
$\alpha$ ranges over multi-indices with 
$|\alpha|=2$, and $\beta$ ranges over multi-indices with $|\beta|=3$)
$$
\frac{\partial^2 f^2}{\partial z^\alpha}\,,\,  \re
\frac{\partial^2
f^1}{\partial z_1\partial w}\,,\,
\frac{\partial^3 f^1}{\partial z^\beta}.\tag3.11
$$
It is straightforward (using ideas similar to those used in the proof
of Theorem B in \S9), and left to the reader, to verify that any
formal mapping \thetag{3.1} that preserves the form (A.i.1) or (A.i.3) of $(M,0)$
can be factored uniquely according to \thetag{3.7} with $T$ and $P$ as
above. We shall say that a choice of $P$, as described above, is a choice of
{\it normalization} for the transformations that preserve the forms
(A.i.1) and (A.i.3), respectively, and that a formal mapping preserving the form
has this normalization if it is factored according to \thetag{3.7}
with this $P$.

Similarly, if $(M,0)$ is of the form (A.i.2), then we shall require
that the polynomials $q_1$, $q_2$ in
\thetag{3.8} are of the form
$$
\aligned
q_1(z,w) =& B_1z_1w+B_2z_2w+\sum_{|\beta|=3} C_\beta z^\beta
+\sum_{|\alpha|=2}D_\alpha z^\alpha w+\\&E_1z_1w^2+E_2z_2w^2+
\sum_{|\beta|=3} F_\beta z^\beta w+Rz_1w^3\\
q_2(z,w) =& G_1z_1w+G_2z_2w+H_1z_1w^2+H_2z_2w^2,\endaligned\tag3.12
$$
for $B_k,C_\alpha,D_\beta,E_k,F_\beta,G_k,H_k\in\bC$ and $R\in\bR$,
and that the formal 
series $f^1$, $f^2$ are such that the constant terms in the following
formal series vanish (the indices $j,k$
below range over $\{1,2\}$, $\alpha$ range over multi-indices with
$|\alpha|=2$, and $\beta$ range over multi-indices with $|\beta|=3$)
$$
\left\{\aligned
&\frac{\partial^2
f^j}{\partial z_k\partial w}\,,\,\frac{\partial^3 f^1}{\partial z^\beta}\,,\,\frac{\partial^3
f^1}{\partial z^\alpha \partial w}\,,\,\frac{\partial^3 f^j}{\partial
z_k\partial w^2}\\&
\frac{\partial^4 f^1}{\partial z^\beta\partial w}\,,\,\re
\frac{\partial^4
f^1}{\partial z_1\partial w^3}. 
\endaligned\right.\tag 3.13
$$
Any transformation preserving the form (A.i.2) can be
factored uniquely according to \thetag{3.7} into such a $P$ and such a
$T$. We say that a choice of $P$, as described above, is a choice of
{\it normalization} for the transformations preserving (A.i.2). 

Now, suppose $(M,0)$ is of one of the forms (A.i.1--3). We write the
equation of $M$ near 0 as follows
$$
\im w=|z_1|^2+p_3(z,\bar z)+F(z,\bar z,\re w).\tag3.14
$$
Here, $p_3(z,\bar z)$ is the homogeneous polynomial of degree 3
corresponding to the form (A.i.1--3), and $F(z,\bar z,\re w)$ is a
real-valued, real-analytic function that is $O(4)$. In what follows,
we shall consider $F(z,\bar z,s)$ as a formal power series 
$$
F(z,\bar z,s)=\sum_{\alpha,\beta,k}c^k_{\alpha\beta} z^\alpha\bar
z^\beta s^k\tag3.15
$$
consisting only of terms of weighted degree greater than 3 (here, $s$
is assigned the weight two) and subjected to the reality condition
$$
c^k_{\alpha\beta}=\overline{ c^k_{\beta\alpha}}.\tag3.16
$$
We shall denote by $\Cal F$ the space of all such power series. 
In order to describe the normal form we shall need to decompose such a
power series $F(z,\bar z,s)$ according to type
$$
F(z,\bar z, s)=\sum_{k,l}F_{kl}(z,\bar z,s),\tag3.17
$$
where $F_{kl}(z,\bar z,s)$ is of type $(k,l)$ i.e. for each $t_1,t_2>0$
$$
F_{kl}(t_1z,t_2\bar z,s)=t_1^kt_2^lF_{kl}(z,\bar z,s).\tag3.18
$$
In what follows, $F_{kl},H_{kl}$, and $N_{kl}$ denote formal
power series of type $(k,l)$. We define the space of normal forms
$\Cal N^1\subset\Cal F$, $\Cal 
N^2\subset\Cal F$, and $\Cal N^3\subset\Cal F$ for
types (A.i.1), (A.i.2) and (A.i.3), respectively, as follows: First,
$N(z,\bar z,s)$ is in {\it regular form} (see the next section for
further discussion of this notion) which can be expressed by 
$$
N(z,\bar z,s)=\sum_{\min(k,l)\geq 1}N_{kl}(z,\bar z,s).\tag3.19
$$
Moreover, the non-zero terms $N_{kl}$ satisfy the following
conditions. In the cases (A.i.1) and (A.i.3), for $j=1,3$,
$$
\aligned
N_{22}\in\Cal N^j_{22}\quad,&\quad N_{32}\in \Cal N^j_{32}\\
N_{42}\in\Cal N^j_{42}\quad,&\quad N_{33}\in \Cal N^j_{33}\\
N_{k1}\in\Cal N^j_{k1}\quad,&\quad k=1,2,3\ldots,\endaligned
\tag3.20
$$
where, for both $j=1,3$,
$$
\aligned
N^j_{11} &= \left\{F_{11}\:F_{11}=z_2H_{01}+\overline{z_2H_{01}}\right\}\\
N^j_{22} &=
\left\{F_{22}\:F_{22}=z^2_2H_{02}+\overline{z^2_2H_{02}}+z_1z_2\bar
z_1\bar z_2 H_{00}\right\}\\
N^j_{33} &= \left\{F_{33}\:F_{33}=z_2H_{23}+\overline{z_2H_{23}}\right\}\\
N^j_{21} &= \left\{F_{21}\:F_{21}=\bar z_2H_{20}\right\},\endaligned\tag3.21
$$
and furthermore
$$
\aligned
N^1_{31} &=
\left\{F_{31}\:F_{31}=\bar z_2 H_{30}+\bar z_1z_1^2H_{10}+\bar z_1z^3_2H_{00}\right\}\\
N^3_{31} &=
\left\{F_{31}\:F_{31}=\bar z_2 H_{30}+\bar z_1z_1^3H_{00}+\bar z_1z^2_2H_{10}
\right\}
,\endaligned\tag3.22
$$
$$
\aligned
N^1_{32} &=
\left\{F_{32}\:F_{32}=z_1^2z_2H^1_{02}+z_2^3H^2_{02}+z_1^3\bar z_2
H_{01}+ z_1z_2^2\bar z_2^2H_{00}\right\}\\
N^3_{32} &=
\left\{F_{32}\:F_{32}=z_1z^2_2H^1_{02}+z_2^3H^2_{02}+z_1^3\bar z_2
H_{01}+z^2_1z_2\bar z^2_2H_{00}\right\}
,\endaligned\tag3.23
$$
$$
\aligned
N^1_{42} &=
\left\{F_{42}\:F_{42}=\bar z^2_2 H_{40}+\bar z^2_1z_2^4H^1_{00}+\bar
z_1\bar z_2(z_1^4H^2_{00}+z_2^4H^3_{00})\right\}\\
N^3_{42} &=
\left\{F_{42}\:F_{42}=\bar z^2_2 H_{40}+\bar z_1^2z_2^4H_{00}+\bar
z_1\bar z_2z^3_2H_{10}
\right\}
,\endaligned\tag3.24
$$
and finally, for $k\geq4$,
$$
\aligned
N^1_{k1} &=
\left\{F_{k1}\:F_{k1}=\bar z_2z_1^{k}H_{00}\right\}\\
N^3_{k1} &=
\left\{F_{k1}\:F_{k1}=\bar z_2z_2^{k}H_{00}
\right\}.\endaligned\tag3.25
$$
If $(M,0)$ instead is of the form (A.i.2), then the terms $N_{kl}$
satisfy the following
$$
\aligned
N_{33}\in\Cal N^2_{33}\quad,&\quad N_{43}\in \Cal N^2_{43}\\
N_{53}\in\Cal N^2_{53}\quad,&\quad N_{44}\in \Cal N^2_{44}\\
N_{54}\in\Cal N^2_{54}\quad,&\quad N_{55}\in \Cal N^2_{55}\\
N_{k1}\in\Cal N^2_{k1}\quad,&\quad N_{k2}\in\Cal N^2_{k2},\quad k=1,2,3\ldots,\endaligned
\tag3.26
$$
where
$$
\aligned
\Cal N^2_{11}
&=\left\{F_{11}\:F_{11}=z_2H_{01}+\overline{z_2H_{01}}\right\}\\
\Cal N^2_{21} &=\left\{F_{21}\:F_{21}=\bar z_2H_{20}+\bar z_1z_2H_{10}\right\}\\
\Cal N^2_{31} &=\left\{F_{31}\:F_{31}=z_2H_{21}+z_1^3\bar
z_2H_{00}\right\}\\
\Cal N^2_{22} &=\left\{F_{22}\:F_{22}=z_2\bar z_2H_{11}\right\}\\
\Cal N^2_{33} &=\left\{F_{33}\:F_{33}=z_2H_{23}+\overline{z_2H_{23}}\right\}\\
\Cal N^2_{43} &=\left\{F_{43}\:F_{43}=\bar z_2H_{42}+\bar z^3_1z^3_2H_{10}\right\}\\
\Cal N^2_{53} &=\left\{F_{53}\:F_{53}=\bar z_2H_{52}+\bar z^3_1z^4_2H_{10}\right\}\\
\Cal N^2_{44}
&=\left\{F_{44}\:F_{44}=z^2_2H_{24}+\overline{z^2_2H_{24}}+z_1^3z_2\bar
z_1^3\bar z_2H_{00}\right\}\\
\Cal N^2_{54} &=\left\{F_{54}\:F_{54}=z_2H_{44}+z^5_1\bar z^2_2H_{02}\right\}\\
\Cal N^2_{55} &=\left\{F_{55}\:F_{55}=z_2H_{45}+\overline{z_2H_{45}}\right\}\\
\Cal N^2_{k1} &=\left\{F_{k1}\:F_{k1}=\bar z_2H_{k0}\right\},\quad k=4,5,\ldots\\
\Cal N^2_{k2} &=\left\{F_{k2}\:F_{k2}=\bar
z_2H_{k1}\right\},\quad k=3,4,\ldots.
\endaligned\tag3.27
$$
We are now in a position to state the theorem on normal forms for
(A.i.1--3).
\proclaim{Theorem B} Let $M$ be a real-analytic hypersurface in
$\bC^3$ given
near $0\in M$ in one of the forms \thetag{A.i.$k$}, for $k\in\{1,2,3\}$, as defined
in Theorem A. Then, given any choice of
normalization (i.e. a choice of $P$ as described above), there is a
unique formal 
transformation \thetag{3.1} with this normalization that transforms the defining equation
\thetag{3.14} of $(M,0)$ to
$$
\im w'=|z'_1|^2+p_3(z',\bar z')+N(z',\bar z',\re w'),\tag3.28
$$
where $N(z,\bar z,s)\in\Cal N^k$. \endproclaim

\flushpar
{\it Remark.} As already mentioned in the remark following Theorem
2.4 above, this result implies an improved bound on the dimension of the
stability group $\Aut(M,p_0)$ of a 2-nondegenerate hypersurface $M$ in
$\bC^3$ of either of the types (A.i.1--3). Counting the number of parameters in the
normalizations, we find that in the case (A.i.1) with $\gamma=0$ we
have $\dim_\bR\Aut(M,p_0)\leq 17$, in the case (A.i.1) with
$\gamma=1$ as well as in the case (A.i.3) we have
$\dim_\bR\Aut(M,p_0)\leq 19$, and in 
the case (A.i.2) we have $\dim_\bR\Aut(M,p_0)\leq 45$.\medskip

The proof of Theorem B will be given in \S9. 
Combining Theorem B with Theorem 2.6 ([BER3]) and Theorem A, we obtain the
following.

\proclaim{Corollary 3.24} Let $M$ and $M'$ be 
real-analytic hypersurfaces in $\bC^3$ that are $2$-nondegenerate and
whose Levi forms have one non-zero eigenvalue at
$p_0\in M$ and $p_0'\in M'$, respectively. Then $(M,p_0)$ and
$(M',p_0')$ are biholomorphically equivalent if and only if, for two
(possibly different) choices of normalization as described in Theorem
B, $(M,p_0)$ and $(M',p_0')$ can be brought to the same normal
form.\endproclaim 


\heading 4. Some applications and examples\endheading
Before turning to the proofs of Theorems A and B, we shall discuss
some applications of these theorems. Therefore, in this section, we
shall assume that the results of \S3 have been proved. 

\subhead 4.1. Computation of two fourth order invariants\endsubhead Let
$M\subset \bC^3$ be a real-analytic hypersurface that is
2-nondegenerate at $p_0\in M$. We shall assume that $(M,p_0)$ is of
the type (A.i.2) (see Theorem A) and compute two fourth order invariants,
corresponding to the Taylor coefficients $c_{22}$ and $b_{22}$ in
\thetag{4.1.2} below, for $M$ at
$p_0$. In fact, the invariants will be
formal invariants
(and then, of course, also biholomorphic invariants) so, in what follows, we
shall work with formal power series and formal transformations. 
In view of Theorems A and B, we may assume that the defining equation
of $M$ is  
$$
\im w=|z_1|^2+z_1^2\bar z_2+\bar
z_1^2z_2+N(z,\bar z,\re w),
\tag 4.1.1
$$
where $N\in\Cal N^2$. We write this in the form
$$
\aligned
\im w = &|z_1|^2+z_1^2\bar z_2+\bar
z_1^2z_2+N_{31}(z,\bar z)+\overline{N_{31}(z,\bar
z)}+\\&|z_2|^2(a_{22}|z_1|^2+b_{22}z_1\bar z_2+\bar b_{22}\bar 
z_1z_2+c_{22}|z_2|^2)+\ldots,\endaligned
\tag 4.1.2
$$
where $a_{22},c_{22}\in\bR$, $b_{22}\in \bC$, $N_{31}\in\Cal
N^2_{31}$, and where the dots $\ldots$ represent terms in the power 
series that are either $O(5)$ or that are $O(4)$ and divisible by
$\re w$. We subject $M$ to a formal invertible transformation
\thetag{3.1} that preserves the normal form, i.e. that corresponds to
another normalization of the transformation to normal form.  Thus, the
transformed formal hypersurface $M'$ is defined by an equation of the
form \thetag{4.1.1} that can be written in the form \thetag{4.1.2}, in
which the corresponding constants and functions 
are denoted with  $'$s. We obtain the following equation for $M'$ in
terms of the transformation $(f^1,f^2,g)$
$$
\aligned
g-\gb= &2i(f^1\fob+(f^1)^2\ftb+\ftb f^2)+\\&2if^2\ftb(a_{22}f^1\fob+b_{22}f^1\ftb+\bar
b_{22} \fob f^2+c_{22}f^2\ftb)+\ldots,\endaligned\tag4.1.3
$$ 
where the dots $\ldots$ this time signify terms that will not influence the constants
$c'_{22}$ and $b'_{22}$. In order to compute $c'_{22}$ and $b'_{22}$,
we shall set $w'=0$ and $\bar w'=\Qb'(\bar z',z,0)$ in \thetag{4.1.3},
where $\bar 
w=\Qb'(\bar z',z',w')$ is the complex defining equation of $M'$ as
defined in \S5 below. By using Proposition 5.7 below and the
form of $\Cal N^2_{11}$, it follows that the coefficients of
$|z_2|^4$, $|z_2|^2z_1\bar z_2$, and $|z_2|^2z_1\bar z_2$ in
$\Qb'(\bar z',z',0)$ are $-2ic_{22}'$, $-2ib'_{22}$, and $-2i\bar
b'_{22}$, respectively. By identifying coefficients on both sides of
\thetag{4.1.3}, in which we have set $w'=0$ and $\bar w'=\Qb'(\bar
z',z',0)$, we obtain
$$
Cc'_{22}=c_{22},\quad
Cb'_{22}=e^{-2i\theta}(2Ac_{22}+C^{1/2}e^{i\theta}b_{22}),\tag 4.1.4
$$
where $C>0$, $\theta\in\bR$, and $A\in D$ are as in Proposition 3.3. (A
similar argument will be needed in the proof of Theorem A and
explained in more detail there; we leave the details of the present
identification to the reader). Defining
$$
\sign(x)=\left\{\aligned \frac{x}{|x|}\,,&\quad x\neq 0\\
0\,\,,&\quad x=0,\endaligned\right.\tag4.1.5
$$
it follows from \thetag{4.1.4} that: 
{\it the integer $\delta_{22}=\sign(c_{22})$ is a formal invariant of
$M$} and, by choosing $C>0$ suitably, we can make $c'_{22}=\delta_{22}$. It
also follows from \thetag{4.1.4} that if
$\delta_{22}\neq 0$, then we can make $b'_{22}=0$. On the other hand, if
$\delta_{22}=0$, then the property $b_{22}\neq0$ is an invariant. We may
therefore define an invariant $\epsilon_{22}$ as follows: $\epsilon_{22}=0$ if
$\delta_{22}\neq 0$ or if $\delta_{22}=0$ and $b_{22}=0$, and
$\epsilon_{22}=1$ if $\delta_{22}=0$ and $b_{22}\neq 0$. By choosing
$A$, $C$, and $\theta$ suitably, 
we can make $b'_{22}=\epsilon_{22}$. Notice that we always have
$\delta_{22}\epsilon_{22}=0$.

\subhead Example 4.1.6\endsubhead The preceding discussion implies
that the following real hypersurfaces are mutually non-equivalent
$$
\aligned
M_1 &=\left\{\im w=|z_1|^2+z_1^2\bar z_2+\bar
z_1^2z_2-|z_2|^4\right\}\\
M_2 &=\left\{\im w=|z_1|^2+z_1^2\bar z_2+\bar z_1^2z_2\right\}\\
M_3 &=\left\{\im w=|z_1|^2+z_1^2\bar z_2+\bar z_1^2z_2+|z_2|^2(z_1\bar
z_2+\bar z_1z_2)\right\}\\
M_4 &=\left\{\im w=|z_1|^2+z_1^2\bar z_2+\bar
z_1^2z_2+|z_2|^4\right\}.\endaligned \tag4.1.7
$$

\subhead 4.2. Everywhere Levi degenerate hypersurfaces\endsubhead
We shall consider real-analytic, everywhere Levi degenerate, 
holomorphically nondegenerate hypersurfaces in $\bC^3$. In view of
Proposition 2.3, such a hypersurface $M$ is 2-nondegenerate
at every $p_0\in M\setminus V$, where $V$ is a proper
real-analytic subset of $M$. 
We assume
that the results in \S3 have been proved and apply these to
$(M,p_0)$ for $p_0\in M\setminus V$. First, however, let us give a couple of
examples of such hypersurfaces. 

\subhead Example 4.2.1. (The light cone)\endsubhead Consider the real hypersurface $M$ defined  
as the set of regular points of 
the cylinder in $\bC^3$ over the light cone in $i\bR^3$, i.e. the set of regular points of
$$
(\im Z_1)^2+(\im Z_2)^2-(\im
Z_3)^2=0.\tag4.2.2
$$
All points $p\in M$ are equivalent via affine transformations. Let
us consider $M$ near the point $p_0=(0,i,i)$. Setting 
$$
Z_1=z_1\quad,\quad Z_2=i+z_2\quad,\quad Z_3=i+w,\tag4.2.3
$$
we find that $M$ is defined near $p_0$, which in the above coordinates
corresponds to $(z_1,z_2,w)=(0,0,0)$, by
$$
w-\bar w=-2i+\sqrt{(z_1-\bar z_1)^2+(z_2-\bar z_2+2i)^2}.\tag4.2.4
$$
Now, Taylor expanding the square root on the right, making a suitable transformation
to regular coordinates, and using the arguments in \S6, we find that 
$(M,p_0)$ (and hence $(M,p)$ for any $p\in M$) can be transformed to a
hypersurface in regular form that has the  
partial normal form (A.i.2); we omit the straightforward
calculations. In particular, it follows that $M$ is  
$2$-nondegenerate (and thus Levi degenerate) at every point.

\subhead Example 4.2.5. (Freeman's hypersurface)\endsubhead M. Freeman gives in [F] the 
following example of a real hypersurface $M$ 
in $\bC^3$ whose Levi form has at least one zero eigenvalue at every
point: $M$ is the set of regular 
points of the real cubic defined by the equation
$$
(\im Z_1)^3+(\im Z_2)^3-(\im
Z_3)^3=0.\tag4.2.6
$$
As remarked in [S], the hypersurface $M$ is holomorphically
nondegenerate and hence, in view of Proposition 2.3, $M$ is
2-nondegenerate on a dense open subset. Since \thetag{4.2.6} is
independent of $\re Z$, the germ $(M,q)$, for any
$q\in M$, is biholomorphically equivalent to $(M,p)$, where $p\in M$
is of the form $p=(iX_1,iX_2,iX_3)$ such that $X_j\in\bR$ and
$$
X_1^3+X_2^3-X_3^3=0.\tag4.2.7
$$
A straightforward calculation shows that if $X_j=0$, for any
$j\in\{1,2,3\}$, then in fact both eigenvalues of the Levi form are
zero. A closer analysis of the defining equation of $M$ shows that $M$
is 3-nondegenerate at such a point. At a point $p\in M$
where $X_1X_2X_3\neq 0$, the calculations are more involved. However,
using the software package {\sl Maple} for some of the symbolic
manipulations, we were able to verify that $M$ at such a point is
2-nondegenerate with exactly one non-zero eigenvalue for the Levi
form and, by following the proof of Theorem A (i) or applying Theorem
4.2.8 below, that $M$ is of the
type (A.i.2) at $p$. 
\medskip

The two examples above both have the same partial normal form,
(A.i.2), at 2-nondegenerate points. The reason for that is the following.

\proclaim{Theorem 4.2.8} Let $M\subset\bC^3$ be a real-analytic
hypersurface and $p_0\in M$. If $M$ is Levi degenerate in an open
neighborhood of $p_0$ and $2$-nondegenerate at $p_0$, then $(M,p_0)$
is biholomorphically equivalent to $(M',0)$, where $M'$ is of the form
\thetag{A.i.2} as defined in Theorem A. \endproclaim

\demo{Proof} We first prove Theorem 4.2.8 under the additional
assumption that the Levi form of $M$ at $p_0$ has
one non-zero eigenvalue. Under this assumption, by Theorem
A, $(M,p_0)$ is equivalent to $(M',0)$, where $M'$ is of one of the
forms (A.i.1--3). We assign the weight one to $z$ and $\bar z$, and
the weight two to $w$ and $\bar w$. We write the defining equation
$\rho'=0$ of $M'$ as follows
$$
\rho'(z,w,\bar z,\bar w)=w-\bar w-2i(|z_1|^2+p_3(z,\bar z))+O(4)=0,\tag4.2.9
$$
where $p_3(z,\bar z)$ is a homogeneous polynomial of degree 3
(determined by the form (A.i.1--3) of $M'$), and
$O(4)$ denotes terms of weighted degree 4 and higher. It is easy to
verify that there is a basis 
$L_1$, $L_2$ for the CR vector fields on $M'$ near $0$ of the form
$$
\aligned
L_1 &=\frac{\partial}{\partial \bar z_1}-2iz_1\frac{\partial}{\partial
\bar w}+r_1(z,\bar z,w,\bar w)\frac{\partial}{\partial \bar w}
\\
L_2 &=\frac{\partial}{\partial \bar z_2}+r_2(z,\bar z,w,\bar
w)\frac{\partial}{\partial \bar w},\endaligned\tag4.2.10
$$
where $r_1$ and $r_2$ are $O(2)$ (cf. also \thetag{4.2.14} below). We obtain (with $Z=(z,w)$)
$$
\aligned
\rho'_Z(Z,\bar Z)= &(-4i\bar
z_1-2ip_{3,z_1}(z,\bar z)+O(3),-2ip_{3,z_2}(z,\bar z)+O(3),1+O(2))\\
L_1\rho'_Z(Z,\bar Z) &=(-4i-2ip_{3,z_1\bar z_1}(z,\bar
z)+O(2),-2ip_{3,z_2\bar z_1}(z,\bar z)+O(2),O(1))\\
L_2\rho'_Z(Z,\bar Z) &=(-2ip_{3,z_1\bar z_2}(z,\bar
z)+O(2),-2ip_{3,z_2\bar z_2}(z,\bar z)+O(2),O(1)),\endaligned\tag4.2.11
$$
where we use the notation
$$
p_{3,z_1}(z,\bar z)=\frac{\partial p_3}{\partial z_1}(z,\bar
z),\quad 
p_{3,z_1\bar z_2}(z,\bar z)=\frac{\partial^2 p_3}{\partial
z_1\partial \bar z_2}(z,\bar z),\quad\text{\rm etc.}\tag4.2.12
$$
The hypersurface $M'$ is Levi degenerate in an open neighborhood
$\omega\subset M'$ of
$0$ if and only if the
three vectors in \thetag{4.2.11} do not span $\bC^3$ for
$(z,w)\in\omega$. Calculating the determinant $D(z,w,\bar z,\bar w)$, we obtain
$$
D(z,w,\bar z,\bar w)=-8p_{3,z_2\bar z_2}(z,\bar z)+O(2).\tag4.2.13
$$
From \thetag{4.2.13} it follows that a necessary condition for
$M'$ to be Levi degenerate in an open neighborhood $\omega$ of $0$ is
that
$$
p_{3,z_2\bar z_2}(z,\bar z)\equiv 0.\tag 4.2.14
$$
A direct calculation shows that $M'$, in this case, must be of the
form (A.i.2). This 
completes the proof under the additional assumption that the Levi form
of $M$ at $p_0$ has a non-zero eigenvalue. 

To complete the proof of
the theorem, in view of Theorem A, we need only to verify that a
real-analytic hypersurface $M'$ of one of the forms (A.ii.1--5) cannot
be Levi degenerate on an open neighborhood of $0$. To do this, we assign the
weight one to $z$ and $\bar z$, and the weight three to $w$ and $\bar
w$. We write the
defining equation of such a hypersurface as follows
$$
\rho'(z,w,\bar z,\bar w)=w-\bar w-2ip_3(z,\bar z)+O(4)=0,\tag4.2.15
$$
where $p_3(z,\bar z)$ is a homogeneous polynomial of degree 3
(determined by the form (A.ii.1--5) of $M'$), and
$O(4)$ denotes terms of weighted degree 4 and higher. A similar
argument (that we leave to the reader) to the one above shows that a
necessary condititon for $M'$ 
to be Levi degenerate on an open neighborhood $\omega$ of $0$ is that
$$
p_{3,z_1\bar z_1}(z,\bar z)p_{3,z_2\bar z_2}(z,\bar z)-
|p_{3,z_1\bar z_2}(z,\bar z)|^2\equiv 0.\tag4.2.16
$$
A direct calculation shows that none of the polynomials  $p_3(z,\bar
z)$ in (A.ii.1--5) satifies \thetag{4.2.16}. This completes the proof.\qed
\enddemo

By including the next term in the weighted homogeneous expansion of
the defining function
$\rho'(z,w,\bar z,\bar w)$ in the calculation above, we would notice
that the invariants $\delta_{22}$ and $\epsilon_{22}$, introduced in
\S4.1 above, at a 2-nondegenerate point of an
everywhere Levi degenerate hypersurface are
necessarily both zero. Indeed, if we denote the term of weighted degree 4 by
$F_4(z,\bar z,s)$, then such a hypersurface must have
$$
F_{4,z_2\bar z_2}(z,\bar z,s)\equiv 4|z_1|^2.\tag4.2.17
$$
Thus, none of the invariants we have computed explicitly can
differ for such hypersurfaces. 

For the proofs of Theorems A and B, we need some
preliminaries. 

\heading 5. Regular coordinates and preliminaries\endheading

We shall say that $(z,w)=(z_1,z_2,w)$ are {\it regular coordinates}
for $(M,p_0)$, where $M$ is a real-analytic hypersurface in $\bC^3$
and $p_0\in M$, if these coordinates vanish at $p_0$ and $M$ can be
described near $p_0=0$  by an equation
$$
\im w=\phi(z,\bar z,\re w),\tag 5.1
$$
where $\phi(z,\chi,u)$ is a holomorphic function near $(0,0,0)$, real-valued for
$(z,\chi,u)=(z,\bar z,s)$ with $s\in\bR$, such that 
$$
\phi(0,\chi,u)\equiv\phi(z,0,u)\equiv0.\tag5.2
$$
Such coordinates always exist (see e.g. [BJT] or [CM]). 
The equation \thetag{5.1} and the hypersurface $M$ are said to be in
{\it regular form} if \thetag{5.2} is satisfied. In the literature
regular form and regular coordinates are sometimes called normal form
and normal coordinates, but in this paper we wish to reserve the 
latter terms for special choices of regular form and regular
coordinates satisfying additional conditions as explained above. 

By writing $\im w=(w-\bar w)/2i$, $\re w=(w+\bar w)/2$, and solving
for $w$ in \thetag{5.1} using the implicit function theorem, we can
describe the hypersurface $M$ near $p_0=0$ by the complex equation
$$
w=Q(z,\bar z, \bar w),\tag5.3
$$
where $Q(z,\chi,\tau)$ is a holomorphic function near $(0,0,0)$. By
complex conjugating \thetag{5.3}, we find that we may also describe
$M$ by
$$
\bar w=\Qb(\bar z,z,w),\tag5.4
$$
where we use the notation $\overline{h}(z)=\overline{h(\bar z)}$. The fact
that the equations \thetag{5.3} and \thetag{5.4} describe a real
hypersurface can be expressed by
$$
\Qb(\chi,z,Q(z,\chi,\tau))\equiv\tau.\tag5.5
$$
Furthermore, equation \thetag{5.2} implies the following 
$$
Q(0,\chi,\tau)\equiv Q(z,0,\tau)\equiv\tau\quad,\quad
\Qb(0,z,w)\equiv\Qb(\chi,0,w)\equiv w.\tag5.6 
$$
We shall refer to equations of the form \thetag{5.3} and \thetag{5.4},
where \thetag{5.6} holds, as being in {\it complex regular form}.
The following proposition is useful.
\proclaim{Proposition 5.7} Let $M$ be given in regular form by
\thetag{5.1} and let \thetag{5.4} be a complex defining equation of
$M$ in complex regular form. Assume that $\phi(z,\bar z,0)$ is
$O(|z|^m)$. (Since $\phi(z,\bar z,s)$ is in regular form, we have 
$m\geq 2$.) If we write $\phi(z,\bar z,s)$ as
$$
\phi(z,\bar z,s)=\phi(z,\bar z,0)+s\hat\phi(z,\bar z,s),\tag 5.8
$$
then 
$$
\Qb(z,\bar z,0)\equiv -2i\,\frac{\phi(z,\bar z,0)}{1+i\hat\phi(z,\bar
z,0)}+O(|z|^{2m+2}).\tag 5.9
$$\endproclaim

\demo{Proof} The function $\Qb(\bar z,z,w)$ is obtained by solving for
$\bar w$ in the equation \thetag{5.1}. If we substitute \thetag{5.4}
in \thetag{5.1} and set $w=0$ (here and in what follows, we consider
$(z,w,\bar z,\bar w)$ as independent variables), then we obtain
$$
\Qb(z,\bar z,0)=-2i\phi(z,\bar z,0)-i\Qb(z,\bar z,0)(\hat\phi(z,\bar
z,0)+O(|z|^2\Qb(z,\bar z,0))).\tag5.10
$$
In the last term on the right, we have used the fact that
$\phi(z,\bar z,s)$ is in regular form to get the factor $|z|^2$. 
Since $\phi(z,\bar z,0)$ is $O(|z|^m)$, we may solve for $\Qb(z,\bar
z,0)$ mod $|z|^{2m+2}$ in \thetag{5.10}. The result is
\thetag{5.9}.\qed\enddemo 

We have the following corollary that suffices for our purposes.
\proclaim{Corollary 5.11} Let $M$, $\phi$, $\Qb$, and $m$ be as in
Proposition $5.7$. Then, we have
$$
\Qb(z,\bar z,0)=-2i\phi(z,\bar z,0)+O(|z|^{m+2}).\tag5.12
$$
\endproclaim
\demo{Proof} If $\hat\phi(z,\bar z,0)$ is $O(|z|^n)$ then Proposition
5.7 implies
$$
\Qb(z,\bar z,0)=-2i\phi(z,\bar z,0)+O(|z|^{p}),\tag5.13
$$
where $p=\min(m+n,2m+2)$. Since $\phi(z,\bar z,s)$ is in regular form,
$n\geq 2$ and Corollary 5.11 follows.\qed\enddemo

In regular coordinates, the following vector fields constitute a basis
for the CR vector fields on $M$ near 0,
$$
L_j=\frac{\partial}{\partial \bar
z_j}+\Qb_{\bar z_j}(\bar z,z,w)\frac{\partial}{\partial \bar w}\quad,\quad
j=1,2,\tag5.14
$$
where, as above, we use the notation
$$
\Qb_{\bar z_j}(\bar z,z,w)=\frac{\partial\Qb}{\partial
\bar z_j}(\bar z,z,w).\tag5.15
$$
It is easy to check that $M$ is $2$-nondegenerate at $p_0=0$ if and
only if the following $5$ vectors span $\bC^2$
$$
\align
\Qb_{z\bar z_1}(0,0,0) &=\left(\Qb_{z_1\bar z_1}(0,0,0)\,,\,\Qb_{z_2\bar
z_1}(0,0,0)\right)\tag5.16\\
\Qb_{z\bar z_2}(0,0,0) &=\left(\Qb_{z_1\bar z_2}(0,0,0)\,,\,\Qb_{z_2\bar
z_2}(0,0,0)\right)\tag5.17\\
\Qb_{z\bar z_1^2}(0,0,0) &=\left(\Qb_{z_1\bar z_1^2}(0,0,0)\,,\,\Qb_{z_2\bar
z_1^2}(0,0,0)\right)\tag5.18\\
\Qb_{z\bar z_2^2}(0,0,0) &=\left(\Qb_{z_1\bar z_2^2}(0,0,0)\,,\,\Qb_{z_2\bar
z_2^2}(0,0,0)\right)\tag5.19\\
\Qb_{z\bar z_1\bar z_2}(0,0,0) &=\left(\Qb_{z_1\bar z_1\bar z_2}(0,0,0)\,,\,\Qb_{z_2\bar
z_1\bar z_2}(0,0,0)\right),\tag5.20
\endalign
$$
but the vectors \thetag{5.16} and \thetag{5.17} do not; if
\thetag{5.16} and \thetag{5.17} span, then $M$ is Levi nondegenerate
(which is the same as $1$-nondegenerate) at $p_0$. 

\heading 6. Proof of Theorem A ({\rm i})\endheading

Let $M$ be a real-analytic hypersurface in $\bC^3$ which is
$2$-nondegenerate at $p_0\in M$ and whose Levi form at $p_0$ has exactly one
non-zero eigenvalue. The first step is the following simple
observation, whose proof is elementary and left to the reader: We may
assume that 
$(z_1,z_2,w)$ are regular coordinates for $(M,p_0)$ and that the
defining equation \thetag{5.1} of $M$ at $p_0=0$ is of the form
$$
\im w=|z_1|^2+O(|z|^3+|\re w||z|^2).\tag6.1
$$
For the purpose of proving Theorem A, it is more convenient to work with the complex
equation \thetag{5.4}, which we write in the following way
$$
\aligned
\bar w= &w-2iz_1\bar z_1+z_1(a_1\bar z_1^2+2b_1\bar z_1\bar z_2+c_1\bar
z_2^2)+\\&z_2(a_2\bar z_1^2+2b_2\bar z_1\bar z_2+c_2\bar
z_2^2)+\ldots,\endaligned\tag6.2
$$
Here, and in what follows, we shall use the following convention
when writing the complex defining equation of a real-analytic
hypersurface: the dots 
$\ldots$ signify terms of degree at least 2 in the
unconjugated variables $(z_1,z_2)$ and terms of total weight at least 4, 
where we assign the weight 1 to the variables $(z_1,z_2,\bar z_1,\bar z_2)$ 
and the weight 2 to $(w,\bar w)$ as in \S3. 
In view of Corollary 5.11, this allows us to reconstruct
the real defining equation of the hypersurface modulo terms of weight at
least $4$. For example, equation \thetag{6.2} implies that the real
defining equation of $M$ is of the form
$$
\aligned
\im w= &|z_1|^2+|z_1|^2\frac{\bar a_1z_1-a_1\bar z_1}{2i}+
|z_1|^2\frac{2\bar b_1z_2-2b_1\bar z_2}{2i}+
\frac{\bar c_1\bar z_1z_2^2-c_1z_1\bar z_2^2}{2i}+\\&
|z_2|^2\frac{\bar c_2z_2-c_2\bar z_2}{2i}+
|z_2|^2\frac{2\bar b_2z_1-2b_2\bar z_1}{2i}+
\frac{\bar a_2\bar z_2z_1^2-a_2z_2\bar z_1^2}{2i}+\\&O(|z|^4+|\re
w||z|^2).\endaligned\tag6.3
$$
With $a_1,a_2,b_1,b_2,c_1,c_2\in\bC$, this is the most general form of
a real equation of the type \thetag{6.1} in regular form. The 5
vectors \thetag{5.16--5.20} become
$$
\align
\Qb_{z\bar z_1}(0,0,0) &=\left(-2i,0\right)\tag6.4\\
\Qb_{z\bar z_2}(0,0,0) &=\left(0,0\right)\tag6.5\\
\Qb_{z\bar z_1^2}(0,0,0) &=\left(2a_1,2a_2\right)\tag6.6\\
\Qb_{z\bar z_2^2}(0,0,0) &=\left(2b_1,2b_2\right)\tag6.7\\
\Qb_{z\bar z_1\bar z_2}(0,0,0) &=\left(2c_1,2c_2\right).\tag6.8
\endalign
$$
Thus, the fact that $M$ is $2$-nondegenerate is equivalent to 
$(a_2,b_2,c_2)$ being different from $(0,0,0)$. We now make a
biholomorphic transformation preserving the origin,
$$
z=f(z',w')\quad,\quad w=g(z',w'),\tag6.9
$$
such that the new coordinates $(z',w')$ are regular for the
transformed hypersurface $(M',0)$. A complex defining equation for the
transformed hypersurface is given by
$$
g(z',w')=Q(f(z',w'),\fb(\bar z',\bar w'),\gb(\bar z',\bar
w')).\tag6.10
$$
It is straightforward to verify that coordinates $(z',w')$ are regular
for a hypersurface $(M',0)$ if and only if for some defining equation 
(and hence for
all defining equations) $\rho(z',w',\bar z',\bar w')=0$ it holds that
$\rho(z',w',0,w')\equiv0$. Thus, $(z',w')$ are regular for $(M',0)$ if
and only if the relation
$$
g(z',w')\equiv Q(f(z',w'),\fb(0,w'),\gb(0,w'))\tag6.11
$$
holds. In particular then,
$$
g(z',w')=w'\hat g(z',w'),\tag6.12
$$
for some holomorphic function $\hat g(z',w')$.
We shall write
$$
g(z',w')=Cw'+w'\hat h(z',w')=Cw'+h(z',w'),\tag6.13
$$
where $C\neq 0$ is a real number, $h(z',w')=w'\hat h(z',w')$, and $\hat h(0,0)=0$.  We also
write, for $j=1,2$,
$$
\aligned
f_j(z',w')=&A^j_1z'_1+A^j_2z'_2+D^jw+B^j_1z_1^2+2B^j_2z_1z_2+B^j_3z_2^2+\\&
O(|z|^3+|w|(|z|+|w|)),\endaligned\tag6.14
$$
where the determinant
$$
A^1_1A^2_2-A^1_2A^2_1\neq0.\tag6.15
$$
We shall show that we can choose 
$A^j_k,B^j_l,C,D^j$ such that the transformed hypersurface $(M',0)$ is 
of one and only one of the model forms \thetag{A.i.1--3}. This will complete
the proof of Theorem A (i).

For notational brevity, we shall drop the $'$ on the new coordinates.
By substituting \thetag{6.9} in \thetag{6.2}, we find the following
equation for the transformed hypersurface $M'$
$$
\aligned
C\bar w+\overline{h}= &Cw+h-2if_1\fb_1+
f_1(a_1\fb_1^2+2b_1\fb_1\fb_2+c_1
\fb_2^2)
+\\&
f_2(a_2\fb_1^2+2b_2\fb_1\fb_2+c_2
\fb_2^2)+\ldots,\endaligned\tag6.16
$$
where we use the convention $\fb_1=\fb_1(\bar z,\bar w)$, $f_1=f_1(z,w)$,
etc. We also use the notation $\ldots$ as explained above, i.e. for those
terms that have degree at least 2 in the unconjugated variables $(z_1,z_2)$
and terms that have total weight at least 4. The new coordinates are regular
for the transformed hypersurface $M'$. We write the complex equation of $M'$ 
in regular form, first as follows
$$
\bar w=\Qb'(\bar z,z,w),\tag 6.17
$$
and then
$$
\aligned
\bar w= &w-2i(r_1z_1\bar z_1+r_2z_1\bar z_2+\bar r_2\bar z_1z_2+r_3
z_2\bar z_2)+z_1(a'_1\bar z_1^2+2b'_1\bar z_1\bar z_2+c'_1\bar
z_2^2)+\\&z_2(a'_2\bar z_1^2+2b'_2\bar z_1\bar z_2+c'_2\bar
z_2^2)+\ldots,\endaligned\tag6.18
$$
We shall identify the coefficients $r_1,r_2,r_3$, and $a'_1,b'_1,c'_1,a'_2,
b'_2,c'_2$ using 
the equation \thetag{6.16}. In order to do this we need to know the linear
part of $\overline{\hat h}(\bar z,0)$. This part is not arbitrary, in view 
of \thetag{6.11}. 
In fact, it follows from that equation that
$$
\overline{\hat h}(\bar z,0)=-2iD^1(\bar A^1_1\bar z_1+\bar A^1_2\bar z_2)+
O(|z|^2).\tag6.19
$$
If we now set $w=0$, substitute $\bar w=\Qb'(z,\bar z,0)$ in \thetag{6.16}, 
and use \thetag{6.19}, we obtain the following equations for the 
new 
coefficients. The coefficients associated with the quadratic part of 
$\Qb'(\bar z,z,0)$ are given
by
$$
\aligned
Cr_1 &= |A^1_1|^2\\
Cr_2 &= A^1_1\bar A^1_2\\
Cr_3 &=|A^1_2|^2.\endaligned\tag6.20
$$
We want to preserve the quadratic terms and, hence, we need to choose
$$
C>0\quad,\quad|A^1_1|=\sqrt{C}\quad,\quad A^1_2=0.\tag6.21
$$
We then get the following equations for the remaining coefficients
$$
\aligned
Ca'_1 =& A^1_1(a_1(\bar A^1_1)^2+2b_1\bar A^1_1\bar A^2_1+c_1(\bar A^2_1)^2
-2i\bar B^1_1)
+\\& A^2_1(a_2(\bar A^1_1)^2+2b_2\bar A^1_1\bar A^2_1+c_2(\bar
A^2_1)^2)+
4D^1\bar A^1_1
\\
Ca'_2 =& A^2_2(a_2(\bar A^1_1)^2+2b_2\bar A^1_1\bar A^2_1+c_2(\bar A^2_1)^2)
\\
Cb'_1 =& 
A^1_1(b_1\bar A^1_1\bar A^2_2+c_1\bar A^2_1\bar A^2_2
-2i\bar B^1_2)
+ 
A^2_1(b_2\bar A^1_1\bar A^2_2+c_2\bar A^2_1\bar A^2_2)
\\
Cb'_2 =& 
A^2_2(b_2\bar A^1_1\bar A^2_2+c_2\bar A^2_1\bar A^2_2)
\\
Cc'_1 =& A^1_1(c_1(\bar A^2_2)^2
-2i\bar B^1_3)
+ A^2_1(c_2(\bar A^2_2)^2)
\\
Cc'_2 =& A^2_2c_2(\bar A^2_2)^2.
\endaligned
\tag6.22
$$
Since $|A^1_1|\neq 0$, we can make $a'_1=b'_1=c'_1=0$ without restricting 
the $A^j_k$ further by a suitable choice of $B^1_1$, $B^1_2$, and $B^1_3$. 

Now, suppose first that $b_2=c_2=0$. 
Then the new coefficients 
$b'_2$ and $c'_2$ are also 0. It follows that $M$ is not biholomorphically 
equivalent to either of \thetag{A.i.1} or \thetag{A.i.3}.
Since $M$ is 2-nondegenerate, $a_2$ 
must be non-zero. It is then easy to see that, by a suitable choice of 
$A^1_1$, $A^2_2$, and $C$, we can make $a'_2=-2i$. Thus, in this case $M$ is
biholomorphically equivalent to \thetag{A.i.2} and not to either
\thetag{A.i.1} or \thetag{A.i.3} . 

Next, suppose $c_2=0$ and $b_2\neq 0$. The coefficient $c'_2$ remains 0
and $b'_2$ cannot be made 0 by a biholomorphic transformation.
By a suitable choice of $A^1_1$, $A^2_2$, and $C$, we can make $b'_2=-i$ and, 
by choosing $A^2_1$ accordingly, we can make $a'_2=0$. Hence, we can
transform $M$ to the form \thetag{A.i.3}, but not to either of \thetag{A.i.1}
or \thetag{A.i.2}. 

Finally, suppose $c_2\neq0$. Then, we cannot make $c'_2=0$ via a 
biholomorphic transformation. However, 
by choosing $A^2_2$ and $C$ such that
$$
-2iC=A^2_2(\bar A^2_2)^2c_2,\tag6.23
$$
we 
can make $c'_2=-2i$. We then make $b'_2=0$ by choosing $\bar A^2_1=
-b_2\bar A^1_1/c_2$. Substituting this in $a'_2$, we find
$$
Ca'_2=A^2_2\bar A^1_1\frac{a_2c_2-b_2^2}{c_2},\tag6.24
$$
where $A^2_2$ also satisfies \thetag{6.23} and
$|A^1_1|=\sqrt{C}$. Now, either $b_2^2-a_2c_2=0$, in which case
$a'_2=0$, or $b_2^2-a_2c_2\neq0$, in which case we can
choose the argument of $A^1_1$ and $C>0$ such that $a'_2=-2i$. Thus,
$M$ can be transformed to the form \thetag{A.i.1}, where
$\gamma=0,1$ is completely determined by $a_2$, $b_2$, $c_2$, and $M$
cannot be transformed to \thetag{A.i.2} or \thetag{A.i.3}. This
completes the proof of Theorem A (i). 

\heading 7. Proof of Theorem A {\rm (ii)}; the beginning\endheading

\subhead 7.1. The setup\endsubhead Let now $M$ be a real-analytic
hypersurface in $\bC^3$ which is  
$2$-nondegenerate at $p_0\in M$ and whose Levi form at $p_0$ is 0. As
in the previous section, we let $(z_1,z_2,w)$ be regular coordinates
for $M$ at $p_0$. We write the complex defining equation \thetag{5.4}
of $M$ as
follows:
$$
\aligned
\bar w= &w+z_1(a_1\bar z_1^2+2b_1\bar z_1\bar z_2+c_1\bar
z_2^2)+\\&z_2(a_2\bar z_1^2+2b_2\bar z_1\bar z_2+c_2\bar
z_2^2)+\ldots,\endaligned\tag7.1.1
$$
We use here the same convention as in \S 3, namely the dots
$\ldots$ signify terms of degree at least 2 in the
unconjugated variables $(z_1,z_2)$ and terms of total weight at least 4. 
As in \S 3 we assign the weight 1 to the variables 
$(z_1,z_2,\bar z_1,\bar z_2)$, but here we assign the weight 3 to 
$(w,\bar w)$. Equation \thetag{5.1} is equivalent to the real equation
for $M$ at $p_0=0$ being of the form
$$
\aligned
\im w= &|z_1|^2\frac{\bar a_1z_1-a_1\bar z_1}{2i}+
|z_1|^2\frac{2\bar b_1z_2-2b_1\bar z_2}{2i}+
\frac{\bar c_1\bar z_1z_2^2-c_1z_1\bar z_2^2}{2i}+\\&
|z_2|^2\frac{\bar c_2z_2-c_2\bar z_2}{2i}+
|z_2|^2\frac{2\bar b_2z_1-2b_2\bar z_1}{2i}+
\frac{\bar a_2\bar z_2z_1^2-a_2z_2\bar z_1^2}{2i}+\\&O(|z|^4+|\re
w||z|^2).\endaligned\tag7.1.2
$$
The first two vectors \thetag{5.16-5.17} are now 0 and the $2$-nondegeneracy
is expressed by the fact that
$$
\align
\Qb_{z\bar z_1^2}(0,0,0) &=\left(2a_1,2a_2\right)\tag7.1.3\\
\Qb_{z\bar z_2^2}(0,0,0) &=\left(2b_1,2b_2\right)\tag7.1.4\\
\Qb_{z\bar z_1\bar z_2}(0,0,0) &=\left(2c_1,2c_2\right).\tag7.1.5
\endalign
$$
span $\bC^2$. 
We now make a biholomorphic transformation
\thetag{6.9} where $(f,g)$ are as in \S 3, i.e. 
$f=(f_1,f_2)$ is of the form \thetag{6.14} and $g$ is
of the form \thetag{6.13} satisfying \thetag{6.10}. 
In this section though,
the coefficients $D^j$ and $B^j_l$ will not enter into the equations as
we shall se. Also, a simple scaling
argument shows that we may assume that $C=1$ in \thetag{6.13}. 
Subjecting $(M,0)$ to this transformation we obtain a new 
real-analytic hypersurface $(M',0)$ for which $(z_1',z_2',w')$ are regular
coordinates and which is given by the complex equation \thetag{6.17}. 
We drop the $'$ on the coordinates and write this equation as
$$
\aligned
\bar w= &w+z_1(a'_1\bar z_1^2+2b'_1\bar z_1\bar z_2+c'_1\bar
z_2^2)+\\&z_2(a'_2\bar z_1^2+2b'_2\bar z_1\bar z_2+c'_2\bar
z_2^2)+\ldots,\endaligned\tag7.1.6
$$
A straightforward calculation shows that
$$
\aligned
a'_1 =& 
A^1_1(a_1(\bar A^1_1)^2+2b_1\bar A^1_1\bar A^2_1+c_1(\bar A^2_1)^2)
+\\& A^2_1(a_2(\bar A^1_1)^2+2b_2\bar A^1_1\bar A^2_1+c_2(\bar A^2_1)^2)
\\
a'_2 =& 
A^1_2(a_1(\bar A^1_1)^2+2b_1\bar A^1_1\bar A^2_1+c_1(\bar A^2_1)^2)
+\\& A^2_2(a_2(\bar A^1_1)^2+2b_2\bar A^1_1\bar A^2_1+c_2(\bar
A^2_1)^2)
\\
b'_1 =& 
A^1_1
(a_1\bar A^1_1\bar A^1_2+b_1(\bar A^1_1\bar A^2_2+\bar A^2_1
\bar A^1_2)+c_1\bar A^2_1\bar A^2_2)
+\\& A^2_1
(a_2\bar A^1_1\bar A^1_2+b_2(\bar A^1_1\bar A^2_2+\bar A^2_1
\bar A^1_2)+c_2\bar A^2_1\bar A^2_2)
\\
b'_2 =& 
A^1_2
(a_1\bar A^1_1\bar A^1_2+b_1(\bar A^1_1\bar A^2_2+\bar A^2_1
\bar A^1_2)+c_1\bar A^2_1\bar A^2_2)
+\\& A^2_2
(a_2\bar A^1_1\bar A^1_2+b_2(\bar A^1_1\bar A^2_2+\bar A^2_1
\bar A^1_2)+c_2\bar A^2_1\bar A^2_2)
\\
c'_1 =& 
A^1_1(a_1(\bar A^1_2)^2+2b_1\bar A^1_2\bar A^2_2+c_1(\bar A^2_2)^2)
+\\& A^2_1(a_2(\bar A^1_2)^2+2b_2\bar A^1_2\bar A^2_2+c_2(\bar A^2_2)^2)
\\
c'_2 =& 
A^1_2(a_1(\bar A^1_2)^2+2b_1\bar A^1_2\bar A^2_2+c_1(\bar A^2_2)^2)
+\\& A^2_2(a_2(\bar A^1_2)^2+2b_2\bar A^1_2\bar A^2_2+c_2(\bar A^2_2)^2)
\endaligned
\tag7.1.7
$$
Before we start the proof of Theorem A (ii), we shall make a preliminary
reduction. Note that the equation for $b'_2$ can be written
$$
\aligned
b'_2= &\bar A^2_1(
b_1|A^1_2|^2+c_1A^1_2\bar A^2_2+b_2\bar A^1_2A^2_2+
c_2|A^2_2|^2)+
\\&\bar A^1_1(
a_1|A^1_2|^2+b_1A^1_2\bar A^2_2+a_2\bar A^1_2A^2_2
+b_2|A^2_2|^2).\endaligned\tag7.1.8
$$
Clearly, we can always choose $A^1_1,A^1_2,A^2_1,A^2_2$ satisfying \thetag
{6.15} such that $b'_2=0$. We shall assume that we have already subjected
$M$ to such a transformation and, consequently, we shall assume in what 
follows that $b_2=0$. 

When we make the transformation \thetag{6.9}, we will
want to keep $b'_2=0$. Thus, in view of \thetag{7.1.8}, the transformation must 
satisfy
the following relation
$$
\bar A^2_1R(A^1_2,\bar A^1_2,A^2_2,\bar A^2_2)+\bar A^1_1
S(A^1_2,\bar A^1_2,A^2_2,\bar A^2_2)=0,\tag 7.1.9
$$
where we use the notation
$$
\aligned
R(u,\bar u,v,\bar v)=&
b_1|u|^2+c_1u\bar v+
c_2|v|^2\\
S(u,\bar u,v,\bar v)=&
a_1|u|^2+b_1u\bar v+a_2\bar uv.\endaligned\tag7.1.10
$$
If $A^1_2, A^2_2$ are chosen such that 
$$
R(A^1_2,\bar A^1_2,A^2_2,\bar A^2_2)\neq0,\tag7.1.11
$$
then 
$$
\bar A^2_1=-\bar A^1_1\frac{S(A^1_2,\bar A^1_2,A^2_2,\bar A^2_2)}
{R(A^1_2,\bar A^1_2,A^2_2,\bar A^2_2)}.\tag7.1.12
$$
We may substitute this into the expression for $b'_1$ to obtain
$$
b'_1=-\bar A^1_1(A^1_1A^2_2-A^2_1A^1_2)\frac{P_1(\bar A^1_2,\bar A^2_2)}{R(A^1_2,\bar
A^1_2,A^2_2,\bar A^2_2)},\tag7.1.13
$$
where
$$
P_1(\bar A^1_2,\bar A^2_2)=-a_2b_1(\bar A^1_2)^2+(a_1c_2-a_2c_1)\bar
A^1_2\bar A^2_2+b_1c_2(\bar A^2_2)^2.\tag7.1.14
$$
If we assume that $A^2_2\neq0$, then we may write
$$
\bar A^1_2=\zeta\bar A^2_2.\tag7.1.15
$$
We then use the notation
$$
\aligned
r(\zeta,\bar \zeta)=& |A^2_2|^{-2}R(
\bar\zeta A^2_2,\zeta \bar A^2_2,
A^2_2,\bar A^2_2)=
b_1|\zeta|^2+c_1\bar\zeta+
c_2\\
s(\zeta,\bar\zeta)=& |A^2_2|^{-2}S(\bar\zeta A^2_2,\zeta \bar A^2_2,
A^2_2,\bar A^2_2)=
a_1|\zeta|^2+b_1\bar \zeta+a_2\zeta
\\
p_1(\zeta)=& (\bar A^2_2)^{-2}P_1(\zeta \bar A^2_2,\bar
A^2_2)=-a_2b_1\zeta^2+
(a_1c_2-a_2c_1)\zeta+b_1c_2
\endaligned\tag7.1.16
$$
so that the relation \thetag{7.1.9} becomes
$$
\bar A^2_1r(\zeta,\bar\zeta)+\bar A^1_1s(\zeta,\bar \zeta)=0.\tag7.1.17
$$
and, provided $r(\zeta,\bar\zeta)\neq0$, 
$$
b'_1=-\bar A^1_1(\bar
A^2_2)^2(A^1_1A^2_2-A^2_1A^1_2)\frac{
p_1(\zeta)}{|A^2_2|^2r(\zeta,\bar\zeta)}.\tag7.1.18
$$
We introduce the determinants
$$
\aligned
\Delta_{ab} &=a_1b_2-a_2b_1=-a_2b_1
\\
\Delta_{bc} &=b_1c_2-b_2c_1=b_1c_2
\\
\Delta_{ac} &=a_1c_2-a_2c_1.\endaligned
\tag7.1.19
$$
so that the $2$-nondegeneracy of $M$ 
means that at least one of these is non-zero. The proof of Theorem A
(ii) will be divided into different cases.

\subhead 7.2 The case $\Delta_{bc}=\Delta_{ac}=0$\endsubhead In
this case, we must have $\Delta_{ab}=-a_2b_1\neq0$. Since the vectors 
$(a_1,a_2)$ and 
$(b_1,0)$ are linearly dependent, $\Delta_{bc}=\Delta_{ac}=0$ implies
that $c_1=c_2=0$.  Our first claim is the following assertion, which implies
that $M$ is not biholomorphically equivalent to either
of the forms \thetag{A.ii.3--5}. 
\proclaim{Assertion 7.2.1} If $b_2=0$, $\Delta_{bc}=\Delta_{ac}=0$ and 
$\Delta_{ab}=-a_2b_1\neq0$
then, after any biholomorphic transformation \thetag{6.9} satisfying 
\thetag{7.1.9}, $a'_2b'_1$ is non-zero.\endproclaim
\demo{Proof} Note first that $R(A^1_2,\bar A^1_2,A^2_2,\bar A^2_2)=b_1
|A^1_2|^2$. Thus, to satisfy \thetag{7.1.9}, either $A^1_2=0$
or \thetag{7.1.12} holds. 

Suppose first that $A^1_2=0$. Then $b'_1=|A^1_1|^2\bar A^2_2b_1$ and  
$a'_2=(\bar A^1_1)^2A^2_2a_2$. This implies that $b'_1$ and $a'_2$ are
non-zero, since
$A^1_2A^2_2$ must be non-zero to make \thetag{6.9} biholomorphic. 

Next, suppose $A^1_2\neq 0$. Then equation \thetag{7.1.12} holds, and
$b'_1$ is given by \thetag{7.1.13}. Since $A^1_2$ is assumed non-zero
and $P_1(\bar A^1_2,\bar A^2_2)=-a_2b_1(\bar A^1_2)^2$, it follows
that $b'_1\neq0$. The equation for $a'_2$ becomes, 
with the notation
$R=R(A^1_2,\bar A^1_2,A^2_2,\bar A^2_2)$ and $S=
S(A^1_2,\bar A^1_2,A^2_2,\bar A^2_2)$,
$$
a'_2=\frac{(\bar A^1_1)^2}{R}
(A^1_2(a_1R-2b_1S)+a_2A^2_2R).
\tag7.2.2
$$
Since $A^1_1\neq 0$ ($A^1_2$ has $A^1_1$ as a factor), $a'_2$ can only
be 0 if
$$
(a_1A^1_2+a_2A^2_2)R-2b_1A^1_2S=0.\tag7.2.3
$$
We claim this is not possible, because \thetag{7.2.3} implies that
\thetag{6.9} is not biholomorphic. To see this, first observe that 
we cannot achieve \thetag{7.2.3} with $A^2_2=0$ unless $a_1=0$. On the
other hand, 
if $a_1=0$, then $A^2_2=0$ implies $A^2_1=0$ which, in turn, implies
that \thetag{6.9} is not biholomorphic. Thus, we may assume $\bar 
A^1_2=\zeta\bar A^2_2$. The equation \thetag{7.2.3} becomes
$$
(a_1\bar \zeta+a_2)r(\zeta,\bar\zeta)-2b_1\bar\zeta s(\zeta,\bar\zeta)=0.
\tag7.2.4
$$
A straightforward calculation shows that \thetag{7.2.4} is equivalent to
$$
\bar \zeta(a_1|\zeta|^2+a_2\zeta+2b_1\bar\zeta)=0\tag7.2.5
$$
Now, 
$$
\aligned
\bar A^1_1\bar A^2_2-\bar A^2_1\bar A^1_2=&
\bar A^1_1\bar A^2_2\left( 
1+\zeta\frac{s(\zeta,\bar\zeta)}{r(\zeta,\bar\zeta)}\right)
\\=&
\frac{\bar A^1_1\bar A^2_2}{r(\zeta,\bar\zeta)}\left( 
r(\zeta,\bar\zeta)+\zeta s(\zeta,\bar\zeta)\right)
\\=&
\frac{\bar A^1_1\bar A^2_2}{r(\zeta,\bar\zeta)}
(a_1|\zeta|^2+a_2\zeta+2b_1\bar\zeta)\zeta.
\endaligned\tag7.2.6
$$
Thus, if the transformation \thetag{6.9} makes $a'_2=0$, then
$\bar A^1_1\bar A^2_2-\bar A^2_1\bar A^1_2=0$, i.e. \thetag{6.9}
is not biholomorphic.
This completes the proof of the assertion.\qed\enddemo

We shall now show that $M$ is equivalent to either \thetag{A.ii.1}, for
some $r>0$, or
\thetag{A.ii.2}. For that we need to make 
$b'_2=c'_1=c'_2=0$ and $b'_1=-i$.
As mentioned above, $b'_2=0$ if either $A^1_2=0$ or \thetag{7.1.12} holds.
First, let us examine the case where $A^1_2\neq0$ and, hence, where \thetag
{7.1.12} holds. It is straightforward to rule out the possibility $A^2_2=0$,
because we cannot make $c'_1=c'_2=0$ via a biholomorphic 
transformation with $A^2_2=0$.
Let us therefore set $\bar A^1_2=\zeta\bar A^2_2$. We obtain the following
equations for $c'_2$, 
$$
\aligned
c'_2 =& A^2_2(\bar A^2_2)^2(\bar\zeta(a_1\zeta^2+2b_1\zeta)+a_2\zeta^2)\\
=& A^2_2(\bar A^2_2)^2(r(\zeta,\bar\zeta)+\zeta s(\zeta,\bar\zeta));
\endaligned
\tag7.2.7
$$ 
the last equality in \thetag{7.2.7} follows easily from the expressions 
\thetag{7.1.16}
for $r$ and $s$. Thus, if $c'_2=0$ then $r(\zeta,\bar\zeta)+\zeta s(\zeta,
\bar\zeta)=0$, since we have already ruled out $A^2_2=0$. However, we also
have
$$
\bar A^1_1\bar A^2_2-\bar A^1_2\bar A^2_1=\bar A^1_2\bar A^2_2\left(1+\zeta
\frac{s(\zeta,\bar\zeta)}{r(\zeta,\bar\zeta)}\right).\tag7.2.8
$$
We deduce that we cannot make $c'_2=0$ by a biholomorphic transformation
if $A^1_2\neq 0$. 

If we choose $A^1_2=0$, then $b'_2=c'_1=c'_2=0$. Our next step is to make
$b'_1=-i$. With $A^1_2=0$, we have
$$
b'_1=|A^1_1|^2\bar A^2_2b_1,\tag7.2.9
$$
so $b'_1=-i$ is accomplished by
$$
\bar A^2_2=\frac{-i}{b_1|A^1_1|^2}.\tag7.2.10
$$
Substituting this in the expression for $a'_2$, we obtain
$$
a'_2=i\left(\frac{\bar A^1_1}{|A^1_1|}\right)^2\frac{a_2}{\bar b_1},\tag7.2.11
$$
By choosing the argument of $\bar A^1_1$
suitably, we can make
$$
a'_2=-2ri,\tag7.2.12
$$
where 
$$
r=\frac{|a_2|}{2|b_1|}.\tag7.2.13
$$
Thus, the number
$r$ is uniquely determined by $a_2$ and $b_1$. 
We also have
$$
a'_1=\bar A^1_1(A^1_1(a_1\bar A^1_1+2b_1\bar A^2_1)+\bar A^1_1A^2_1a_2),\tag
7.2.14
$$
or with $\bar A^2_1=\zeta\bar A^1_1$,
$$
a'_1=|A^1_1|^2\bar A^1_1(a_1+2b_1\zeta+a_2\bar\zeta).\tag7.2.15
$$
We need the following lemma, whose proof is elementary and left to
the reader.
\proclaim{Lemma 7.2.16} The equation $L(\zeta,\bar\zeta)=0$, where
$$
L(\zeta,\bar\zeta)=a_1+2b_1\zeta+a_2\bar\zeta,\tag7.2.17
$$
has a solution if and only if one of the following holds: 
$$
|a_2|\neq 2|b_1|,\tag7.2.18
$$
$$
a_2=2b_1e^{i\theta}\quad\text{\rm and}\quad \frac{a_1}{a_2}
e^{i\theta/2}\in\bR.\tag7.2.19
$$
Moreover, if
$$
a_2=2b_1e^{i\theta},\tag7.2.20 
$$
then the range of $e^{i\theta/2}L(\zeta,\bar\zeta)/a_2$ equals
$$
\left\{x+i\,\im\left(\frac{a_1}{a_2}e^{i\theta/2}\right)\:x\in\bR\right\}.
\tag
7.2.21
$$
\endproclaim
Writing $a_2=2rb_1e^{i\theta}$,  with
$\theta\in[0,2\pi)$, and using that the argument of $A^1_1$
is chosen so as to make $a'_2=-2ri$, we can rewrite \thetag{7.2.15} as follows
$$
a'_1=\pm i|A^1_1|^3|b_1|\frac{e^{i\theta/2}L(\zeta,\bar
\zeta)}{a_2},\tag 7.2.22
$$
where the plus or minus sign corresponds to which branch of the square root was chosen in
defining the argument of $A^1_1$. Now, if the invariant $r$, defined
by \thetag{7.2.13} 
above, does not equal
1, then $M$ is not biholomorphically equivalent to (A.ii.2). On the
other hand, in this case \thetag{7.2.18} holds and, hence,
$L(\zeta,\bar \zeta)$ has a root. Consequently, we can make $a'_1=0$
and $M$ is equivalent to (A.ii.1) with $r$ defined by
\thetag{7.2.13}. 

If the invariant $r=1$, then \thetag{7.2.20} holds. Assume first that
the additional condition in 
\thetag{7.2.19} holds. Then $a'_1$ is always purely
imaginary, in view of Lemma 7.2.16 and \thetag{7.2.22}, and it follows
that $M$ is not equivalent to (A.ii.2). On the other hand, in this
case Lemma 7.2.16 also asserts that $L(\zeta,\bar \zeta)$ has a
root. Thus, we can make $a'_1=0$ and $M$ is equivalent to (A.ii.1)
with $r=1$. 

Finally, if \thetag{7.2.20} holds but not the additional condition
in \thetag{7.2.19}, then we cannot, in view of Lemma 7.2.16 and \thetag{7.2.22}, make
$a'_1$ vanish. Thus, $M$ is not equivalent to (A.ii.1). On the other
hand, by choosing $\zeta$ and the modulus of $A^1_1$ appropriately,
and the
plus or minus sign in \thetag{7.2.22} accordingly, we can make
$a'_1=-2$. It follows that $M$ is equivalent to (A.ii.2). 
This completes the proof of  
Theorem A (ii) in the case $\Delta_{bc}=\Delta_{ac}=0$.

\subhead 7.3. The case $\Delta_{ab}=\Delta_{ac}=0$\endsubhead In
this case, we must have $\Delta_{bc}=b_1c_2\neq0$. As in the previous
case, the fact that $\Delta_{ab}=\Delta_{ac}=0$ implies
that $a_1=a_2=0$. The first step is the following assertion, which shows
that $M$ is not equivalent to \thetag{A.ii.1}, \thetag{A.ii.2},
\thetag{A.ii.4}, or \thetag{A.ii.5}. The proof of this is very similar to the proof of
Assertion 7.2.1 above, and is therefore left to the reader.

\proclaim{Assertion 7.3.1} If $b_2=0$, $\Delta_{ab}=\Delta_{ac}=0$ and 
$\Delta_{bc}=b_1c_2\neq0$
then, after any biholomorphic transformation \thetag{6.9} satisfying 
\thetag{7.1.9}, $b'_1c'_2$ is non-zero.\endproclaim

We are left with showing that $M$ is equivalent to \thetag{A.ii.3}
for precisely one choice of $\lambda\in\bC$, $\lambda\neq0$. 
One possible way of satisfying \thetag{7.1.9} would be to choose
$A^1_2$ and $A^2_2$ such that $R(A^1_2,\bar A^1_2,A^2_2,\bar A^2_2)=0$
and $A^1_1=0$. However, it is easy to see that such a transformation
leads to a hypersurface which is not of any of the forms \thetag{A.ii.1--5}. 
The only alternative in this case
is to choose $\bar A^2_1$ according to
\thetag{7.1.12}. In order for the transformation to be biholomorphic,
we must also choose $\bar A^1_1\neq 0$ and $\bar A^2_2\neq 0$. 
We wish to make a transformation of this form such that 
$a'_1=a'_2=0$. We set $A^1_2=\zeta\bar A^2_2$ and substitute
\thetag{7.1.12} in the equation for $a'_2$. We obtain
$$
\aligned
a'_2=& (\bar A^1_1)^1A^2_2\left(\bar
\zeta\left(-2b_1\frac{s(\zeta,\bar\zeta)}{r(\zeta,\bar\zeta)}+
c_1\left(\frac{s(\zeta,\bar\zeta)}{r(\zeta,\bar\zeta)}\right)^2\right)+c_2
\left(\frac{s(\zeta,\bar\zeta)}{r(\zeta,\bar\zeta)}\right)^2\right)\\
=& (\bar
A^1_1)^2A^2_2\frac{s(\zeta,\bar\zeta)}{(r(\zeta,\bar\zeta))^2}
\left((c_2+\bar\zeta c_1)s(\zeta,\bar\zeta)-2b_1\bar\zeta
r(\zeta,\bar\zeta)\right)\\
=& -2b_1\bar\zeta(\bar
A^1_1)^2A^2_2\frac{s(\zeta,\bar\zeta)}{(r(\zeta,\bar\zeta))^2}
\left(r(\zeta,\bar\zeta)+\zeta s(\zeta,\bar\zeta)\right).\endaligned
\tag7.3.2
$$
Also, as in \thetag{7.2.6}, we have
$$
\bar A^1_1\bar A^2_2-\bar A^2_1\bar A^1_2=\frac{\bar A^1_1\bar
A^2_2}{r(\zeta,\bar\zeta)}(r(\zeta,\bar\zeta)+\zeta s(\zeta,\bar\zeta)).\tag7.3.3
$$
It follows that the only way to obtain $a'_2=0$, via a biholomorphic
transformation of the form considered here, is to choose $\zeta=0$,
i.e. $A^1_2=0$. This in turn 
forces $A^2_1=0$, since $s(0,0)=0$. Moreover, $A^2_1=0$ implies $a'_1=0$.
We are left with the following equations
$$\aligned
b'_1 =& |A^1_1|^2\bar A^2_2 b_1\\c'_1 =& A^1_1(\bar A^2_2)^2c_1\\c'_2
=& |A^2_2|^2\bar A^2_2c_2.\endaligned\tag7.3.4
$$
We may choose $\bar A^1_1$ and $\bar A^2_2$
such that $b'_1=-i$ and $c'_1=-2i$. In doing so, we obtain
$c'_2=-2i\lambda$, where $\lambda=2\bar b_1c_2/|c_1|^2$. Thus, $M$ is
equivalent
to \thetag{A.ii.3}, where $\lambda\in\bC$, $\lambda\neq 0$, is uniquely
determined by $b_1$, $c_1$, and $c_2$. This completes the proof of Theorem A (ii) in the
case
$\Delta_{ab}=\Delta_{ac}=0$.

\subhead 7.4. The case $\Delta_{ab}=\Delta_{bc}=0$\endsubhead In
this case, we have $\Delta_{ac}\neq 0$ and $(b_1,b_2)=(0,0)$. The
following shows that $M$ is not equivalent to any of the forms
\thetag{A.ii.1-3}.
\proclaim{Assertion 7.4.1} If $b_1=b_2=0$ and
$\Delta_{ac}=a_1c_2-a_2c_1\neq0$ then, after any biholomorphic
transformation
\thetag{6.9} satisfying \thetag{7.1.9}, $(a'_1,a'_2)\neq(0,0)$ and 
$(c'_1,c'_2)\neq(0,0)$.\endproclaim

\demo{Proof} We begin by assuming that $a_1c_2\neq 0$. First, let us
set $\bar A^2_2=0$. In 
order to satisfy \thetag{7.1.9} and still have a biholomorphic
transformation, we must make $\bar A^1_1=0$. It is then easy to check that
$a'_1c'_2\neq 0$. 

Next, let us assume $\bar A^2_2\neq0$ and make $R(A^1_2,\bar
A^1_2,A^2_2,\bar A^2_2)=0$ by choosing $A^1_2$,
$A^2_2$ such that
$$
c_1A^1_2+c_2A^2_2=0.\tag7.4.2
$$
In order to satisfy \thetag{7.1.9} and still have a biholomorphic
transformation, we must have $\bar A^1_1=0$, because $\bar A^2_2\neq 0$, 
$a_1c_2-a_2c_1\neq 0$, and \thetag{7.4.2} imply that
$a_1A^1_2+a_2A^2_2 \neq 0$ and $\bar A^1_2\neq 0$. Such a 
transformation yields $a'_1\neq0$, as is easy to 
verify. Solving for $A^2_2$ in \thetag{7.4.2} and substituting in the
expression for $c'_2$, we obtain
$$
c'_2=\frac{A^1_2(\bar A^1_2)^2}{c_2}(a_1c_2-a_2c_1),\tag7.4.3
$$
which cannot be made to vanish when $\bar A^1_1=0$ since the latter implies that
$\bar A^1_2\neq 0$.

Thirdly, let us consider the situation where $R(A^1_2,\bar
A^1_2,A^2_2,\bar A^2_2)\neq 0$ and the
equation \thetag{7.1.9} is  
satisfied by choosing  
$\bar A^2_1$ according to \thetag{7.1.12}. Setting
$\bar A^1_2=\zeta\bar A^2_2$ and substituting \thetag{7.1.12} in the
equation for $a'_2$ we obtain
$$
\aligned
a'_2 =& A^2_2(\bar
A^1_1)^2\left(\bar\zeta\left(a_1+c_1\left(\frac{s(\zeta,\bar\zeta)}
{r(\zeta,\bar\zeta)}\right)^2\right)+a_2+c_2\left(\frac{s(\zeta,\bar\zeta)}
{r(\zeta,\bar\zeta)}\right)^2\right)\\=&
\frac{A^2_2(\bar
A^1_1)^2}{r(\zeta,\bar\zeta)^2}(c_1\bar\zeta+c_2)(a_1\bar\zeta+a_2)(r(\zeta,\bar
\zeta)+\zeta s(\zeta,\bar\zeta)).\endaligned\tag7.4.4
$$
As we have seen above, it is not possible to make
$r(\zeta,\bar\zeta)+\zeta s(\zeta,\bar\zeta)=0$
by a biholomorphic transformation (see the proofs of Assertion 7.2.1
and 7.3.1). Moreover, $r(\zeta,\bar\zeta)$ must be different from
0, which is equivalent to $c_1\bar\zeta+c_2\neq0$. Hence, to make $a'_2=0$ we
must make $a_1\bar\zeta+a_2=0$. This in turn implies that
$s(\zeta,\bar\zeta)=0$.
Substituting this in the equation for $a'_1$, we obtain
$a'_1=A^1_1(\bar A^1_1)^2 a_1$. The latter cannot be made 0
by a biholomorphic transformation. Thus, we cannot make both
$a'_1=0$ and $a'_2=0$. A similar argument shows that we cannot make
both $c'_1=0$ and $c'_2=0$. This concludes the proof of Assertion
7.4.1 under the additional assumption that $a_1c_2\neq 0$. If
$a_1c_2=0$, then $a_2c_1\neq 0$. The proof in this case is similar to
the one above and left to the reader.\qed\enddemo

Now, to transform $M$ into one of the forms \thetag{A.ii.4} or
\thetag{A.ii.5} we need to make $b'_1=0$, in addition to satifying
\thetag{7.1.9}. We claim that this implies that the only possible
biholomorphic transformations have either $A^1_1=A^2_2=0$ or $A^1_2=A^2_1=0$. We
shall prove it in the case $a_1=c_2=0$. The remaining cases are
similar and left to the reader. Thus, assume $a_1=c_2=0$ which, in
turn, implies $a_2c_1\neq0$. Then,
\thetag{7.1.9} can be satisfied in three ways, either $\bar A^2_2=0$,
$\bar A^1_2=0$, or \thetag{7.1.12} holds. If $\bar A^2_2=0$, then we
need to choose $\bar A^1_1=0$ in order to make $b'_1=0$. If $\bar
A^1_2=0$, then we need $\bar A^2_1=0$ to make $b'_1=0$. If we choose
$\bar A^2_1$ according to \thetag{7.1.12}, assuming then that
$R(A^1_2,\bar A^1_2,A^2_2\bar A^2_2)\neq 0$, then $b'_1$ is given by
\thetag{7.1.13}. This cannot be made 0 via a biholomorphic
transformation, as is easy to verify, since $R(A^1_2,\bar
A^1_2,A^2_2\bar A^2_2)\neq 0$ implies $A^1_2A^2_2\neq 0$ and $\bar
A^1_1=0$ would imply $\bar A^2_1=0$ which is not possible for a
biholomorphic transformation. This proves the claim above. 

Now, any transformation with $A^1_1=A^2_2=0$ or
$A^1_2=A^2_1=0$ satisfies \thetag{7.1.9} and makes $b'_1=0$. Moreover,
for a
transformation with $A^1_1=A^2_2=0$ we obtain
$$
\aligned
a'_1=c_2A^2_1(\bar A^2_1)^2\quad &,\quad a'_2=c_1A^1_2(\bar
A^2_1)^2\\
c'_1=a_2A^2_1(\bar A^1_2)^2\quad &,\quad c'_2=a_1A^1_2(\bar
A^1_2)^2,\endaligned
\tag7.4.5
$$
whereas for a transformation with $A^1_2=A^2_1=0$, we obtain
$$
\aligned
a'_1=a_1A^1_1(\bar A^1_1)^2\quad &,\quad a'_2=a_2A^2_2(\bar
A^1_1)^2\\
c'_1=c_1A^1_1(\bar A^2_2)^2\quad &,\quad c'_2=c_2A^2_2(\bar
A^2_2)^2.\endaligned
\tag7.4.6
$$
Note that in both cases $a'_1c'_2$ and $a'_2c'_1$ differ from
$a_1c_2$ and $a_2c_1$, respectively, by a non-zero constant. Thus,
the properties $a_1c_2\neq0$ and $a_2c_1\neq 0$ are invariant. 

Suppose first that $a_2c_1\neq 0$ and $a_1c_2=0$. Since $a_1c_2=0$ and
this property is invariant, $M$ is not equivalent to (A.ii.4) for any
choices of $\mu,\nu\in\bC$. On the other hand, by choosing
the appropriate transformation, either $A^1_1=A^2_2=0$ or
$A^1_2=A^2_1=0$, and the non-zero coefficients suitably, we can make
$c'_2=0$, $a_2=c_1=-2i$. This makes $a_1=-2i\eta$, where $\eta\in \bC$ is
uniquely determined by $a_1$, $a_2$, $c_1$, and $c_2$. This proves
that $M$ is equivalent to (A.ii.5) for precisely one $\eta\in\bC$.

Suppose next that $a_1c_2\neq 0$. Then $M$ is not equivalent to
(A.ii.5) for any $\eta\in\bC$. On the other hand, by choosing
$A^1_2=A^2_1=0$ and $A^1_1$, $A^2_2$ appropriately, we can make
$a'_1=c'_2=-2i$. This makes $a_2=-2i\alpha$ and $c_1=-2i\beta$, where
$\alpha,\beta\in\bC$ are completely determined by $a_1$, $a_2$, $c_1$, and
$c_2$. Since $a'_1c'_2-a'_2c'_1\neq 0$, it is easy to see that
$\alpha\beta\neq 1$.

However, in this case we could also choose the other type of
transformation, $A^1_1=A^2_2=0$. By choosing $A^1_2$, $A^2_1$
appropriately, we can make $a_1=c_2=-2i$. As is easy to verify, this
makes $a_2=-2i\beta$ and $c_1=-2i\alpha$, where $\alpha$ and $\beta$ are the
same numbers as above. Thus, $M$ can be transformed to (A.ii.4) for
precisly two choices of the pair $(\mu,\nu)$, with $\mu\nu\neq 1$, namely
$(\alpha,\beta)$ and $(\beta,\alpha)$. We can make the choice unique
by requiring, e.g., that $|\mu|\geq|\nu|$ and $\arg\mu\geq\arg\nu$,
where $\arg\mu,\arg\nu\in[0,2\pi)$, if 
$|\mu|=|\nu|$. This
completes the proof of Theorem A (ii) in the case
$\Delta_{ab}=\Delta_{bc}=0$.

\heading 8. Proof of Theorem A {\rm (ii)}; the conclusion\endheading

\subhead 8.1. Two simple lemmas\endsubhead We shall keep the notation
established in the preceeding sections. 
We start by stating two simple, but useful lemmas. Their validity
is easy to check and the details are left to the reader.
\proclaim{Lemma 8.1.1} The following identity holds
$$
p_1(\zeta)=(a_1\zeta+b_1)r(\zeta,\bar\zeta)-(b_1\zeta+c_1)s(\zeta,\bar\zeta).\tag8.1.2
$$
\endproclaim

\proclaim{Lemma 8.1.3} If $\zeta_0$ is a root of $p_1(\zeta)$ such that
$r(\zeta_0,\bar\zeta_0)\neq 0$ then, with $\bar A^2_1$ chosen
according to \thetag{7.1.12} and $\bar A^1_2=\zeta_0\bar A^2_2$, the
following
holds
$$
\bar A^1_1\bar A^2_2-\bar A^1_2\bar A^2_1=\frac{\bar A^1_1\bar
A^2_2}{b_1\zeta_0+c_1}(a_1\zeta_0^2+2b_1\zeta_0+c_1).\tag8.1.4
$$
\endproclaim
\flushpar
We now proceed with the proof of Theorem A.
\subhead 8.2. The case $\Delta_{ac}=0$ and
$\Delta_{ab}\Delta_{bc}\neq0$\endsubhead Since $b_2=0$, it follows
that $a_2b_1$ and $b_1c_2$ are both non-zero. We shall make a biholomorphic
transformation such that $b'_1=b'_2=0$. This will complete the proof
in this case, because the transformed hypersurface $M'$ was treated in
section 7.4 above. We set 
$$
\bar A^1_2=\zeta\bar
A^2_2.\tag 8.2.1
$$
Equation \thetag{7.1.9}, which is equivalent to $b'_2=0$, becomes
$$
\bar A^2_1r(\zeta,\bar\zeta)+\bar A^1_1 s(\zeta,\bar\zeta)=0.\tag8.2.2
$$
If $r(\zeta,\bar\zeta)=0$ and $s(\zeta,\bar\zeta)=0$ have a common
solution $\zeta_0$, then \thetag{8.2.2} can be satisfied by choosing
$\zeta=\zeta_0$. If they do not, then we satisfy \thetag{8.2.2} by
$$
\bar A^2_1=-\bar
A^1_1\frac{s(\zeta,\bar\zeta)}{r(\zeta,\bar\zeta)},\tag8.2.3
$$ 
where $\zeta$ is any complex number such that $r(\zeta,\bar\zeta)\neq0$. (Equation
\thetag{8.2.3} is just \thetag{7.1.12} with $\bar A^1_2=\zeta\bar
A^2_2$). Substituting 
\thetag{8.2.3} in $b'_1=0$, we obtain \thetag{7.1.18}. Note that in
this case
$$
p_1(\zeta)=b_1(-a_2\zeta^2+c_2)\tag8.2.4
$$
and, hence, $p_1(\zeta)$ has two distinct roots 
$$
\zeta_1=\sqrt{\frac{c_2}{a_2}}\quad,\quad\zeta_2=-\sqrt{\frac{c_2}{a_2}}.\tag8.2.5
$$ 
Suppose first that one of those roots, say $\zeta_1$ (the other case
being similar), is not a
solution of 
$$
r(\zeta,\bar\zeta)=0.\tag 8.2.6
$$ 
Then we can make $b'_1=b'_2=0$ by choosing $\bar A^1_2$ and $\bar A^2_1$ according to
\thetag{8.2.1} and \thetag{8.2.3}, respectively, and setting $\zeta=\zeta_1$. We just need to
check that this transformation is biholomorphic. In view of Lemma
8.1.3, it suffices to check that
$$
a_1\zeta_1^2+2b_1\zeta_1+c_1\neq0.\tag8.2.7
$$
Using \thetag{8.2.5} and the fact that $a_2c_1=a_1c_2$, we obtain
$$
\aligned
a_1\zeta_1^2+2b_1\zeta_1+c_1 =&\,
c_1\left(\frac{a_2}{c_2}\zeta_1^2+1\right)+2b_1\zeta_1
\\ =&\,
2(c_1+b_1\zeta_1).\endaligned\tag8.2.8
$$
Since $\zeta_1$ is assumed not to be a solution of \thetag{8.2.6}, it
follows from \thetag{8.1.2} that either \thetag{8.2.7} holds, in which
case we are done, or $\zeta_1$ is a root of both $b_1\zeta+c_1$ and
$a_1\zeta+b_1$, i.e. 
$$
\frac{c_1}{b_1}=\frac{b_1}{a_1}\tag8.2.9
$$
or equivalently
$$
b_1^2-a_1c_1=0.\tag8.2.10
$$
If $\zeta_2$ is also not a solution of \thetag{8.2.6} then we
could have chosen $\zeta=\zeta_2$ instead of $\zeta=\zeta_1$ and we
would be done, since at most one of the points $\zeta_1$, $\zeta_2$ can be
a root of $c_1+b_1\zeta$. Thus, we are left with the situation where
\thetag{8.2.10} holds, and where one of the roots $\zeta_1$, $\zeta_2$
is also a solution of \thetag{8.2.6}. The identity \thetag{8.2.10} implies
that
$$
\zeta_1=-\frac{c_1}{b_1}=-\frac{b_1}{a_1}\quad,\quad
\zeta_2=\frac{c_1}{b_1}=\frac{b_1}{a_1}.\tag8.2.11
$$
A simple computation shows that
$$
r(\zeta_1,\bar\zeta_1)=c_2\quad,\quad
s(\zeta_1,\bar\zeta_1)=-\frac{a_2b_1}{a_1},\tag8.2.12
$$
so the root which is also a solution of \thetag{8.2.6} can only be
$\zeta_2$. It follows from Lemma 8.1.1 that
$s(\zeta_2,\bar\zeta_2)=0$. We can then satisfy \thetag{8.2.2} by
choosing $\bar A^1_2=\zeta_2\bar A^2_2$. Setting 
$$
\bar A^1_1=t\bar A^2_1\tag8.2.13
$$
and substituting in the equation for $b'_1$, we obtain
$$
\aligned
b'_1 =& A^2_1\bar A^2_2(\bar
t(a_1\zeta_2+t+b_1(t+\zeta_2)+c_1)+a_2t\zeta_2+c_2)\\
=& A^2_1\bar A^2_2(\zeta_2(a_1|t|^2+b_1\bar t+a_2t)+b_1|t|^2+c_1\bar
t+c_2)\\=& A^2_1\bar A^2_2(r(t,\bar t))+\zeta_2s(t,\bar t)).\endaligned\tag8.2.14
$$
We claim that $b'_1=0$ if $t=-c_1/b_1$. To see this, we substitute
this value for $t$ in \thetag{8.2.14} and use 
\thetag{8.2.11} and \thetag{8.2.12} to find
$$
\aligned
b'_1=& A^2_1\bar
A^2_2\left(c_2-\frac{c_1}{b_1}\,\frac{a_2b_1}{a_1}\right)
\\=&
A^2_1\bar
A^2_2\frac{1}{a_1}\left(a_1c_2-a_2c_1\right),\endaligned\tag8.2.15
$$
and $a_1c_2-a_2c_1=0$ by assumption. 

Next suppose that both roots of $p_1(\zeta)$ solve \thetag{8.2.6}. In
view of \thetag{8.2.12} above, we observe that \thetag{8.2.6} and 
$$
b_1\zeta+c_1=0\tag 8.2.16
$$
have no common solutions. Thus, it follows from Lemma 8.1.1
that both roots of $p_1(\zeta)$ are common solutions of \thetag{8.2.6}
and
$$
s(\zeta,\bar \zeta)=0.\tag8.2.17
$$
We may then satisfy \thetag{8.2.2} by choosing e.g. $\bar
A^1_2=\zeta_1\bar A^2_2$. If we substitute this and
\thetag{8.2.13} in the equation for
$b'_1$, we obtain 
$$
b'_1 = A^2_1\bar A^2_2(r(t,\bar t))+\zeta_1s(t,\bar t)).\tag8.2.18
$$
Choosing $t=\zeta_2$, we obtain $b'_1=0$. We need to check that the
corresponding transformation is biholomorphic. This
is clear because
$$
\bar A^1_1\bar A^2_2-\bar A^1_2\bar A^2_1 = \bar A^2_1\bar
A^2_2(\zeta_2-\zeta_1)\tag8.2.19
$$
and $\zeta_2\neq\zeta_1$. This completes the proof in the case where $\Delta_{ac}=0$ and
$\Delta_{ab}\Delta_{bc}\neq0$.

\subhead 8.3. The case $\Delta_{ab}=0$ and
$\Delta_{ac}\Delta_{bc}\neq0$\endsubhead As above, it suffices to make a biholomorphic
transformation such that $b'_1=b'_2=0$. Now, it follows from the
assumptions that $a_2=0$ and $a_1b_1c_2\neq0$. We leave it to the
reader to verify that a transformation \thetag{6.9} with $\bar
A^2_2=0$ and $\bar A^2_1=-\bar A^1_1a_1/b_1$ makes $b'_1=b'_2=0$. The
parameter $\bar A^1_2$ can be chosen arbitrarily and the
transformation is biholomorphic as long as $\bar A^1_1\bar A^1_2\neq 0$.

\subhead 8.4. The case $\Delta_{bc}=0$ and
$\Delta_{ab}\Delta_{ac}\neq0$\endsubhead We shall make a biholomorphic
transformation such that $b'_1=b'_2=0$. As above, it follows from the
assumptions that $c_2=0$ and $a_2b_1c_1\neq0$. It is easy to verify
that a transformation \thetag{6.9} with $\bar 
A^1_1=0$ and $\bar A^1_2=-\bar A^2_2c_1/b_1$ makes $b'_1=b'_2=0$. Such
a transformation is biholomorphic as long as $\bar A^2_1\bar A^2_2\neq 0$.

\subhead 8.5. Last case; 
$\Delta_{ab}\Delta_{ac}\Delta_{bc}\neq0$\endsubhead We have
$a_2b_1c_2\neq 0$ and at least one of $c_1$, $a_1$ is
non-zero. This time it suffices to make a biholomorphic transformation
\thetag{6.9}, satisfying \thetag{7.1.9}, such that
the transformed hypersurface $M'$ falls into one of the six
categories considered in sections 7.2--7.4 and 8.2--8.4. Thus, it
suffices to make $b'_1=0$, $a'_2=0$, $c'_2=0$, or $(a_1,c_1)=(0,0)$. 

First, let us assume that one root $\zeta_0$ of $p_1(\zeta)$ is also a
solution of \thetag{8.2.6}. Since $c_2\neq 0$, equation
\thetag{8.2.12} and Lemma 8.1.1 imply that $\zeta_0$ is a solution of
\thetag{8.2.17} as well. Thus, \thetag{8.2.6} and \thetag{8.2.17} have
a common solution $\zeta_0$. As above, this means that we may satisfy
\thetag{7.1.9} by choosing $\bar A^1_2=\zeta_0\bar A^2_2$, without making
any choices of $\bar A^1_1$ and $\bar A^2_1$. We then have
$$
\aligned
c'_2 =&\, A^2_2(\bar A^2_2)^2(\bar
\zeta_0(a_1\zeta_0^2+2b_1\zeta_0+c_1)+a_2\zeta_0^2+c_2)\\=&\,
A^2_2(\bar
A^2_2)^2(r(\zeta_0,\bar\zeta_0)+\zeta_0s(\zeta_0,\bar\zeta_0))\\=&\,
0.\endaligned\tag8.5.1 
$$
Consequently, we are left with the situation where no root of
$p_1(\zeta)$ solves \thetag{8.2.6}. Let us denote the roots by
$\zeta_1$, $\zeta_2$, and observe that these two roots need not be
distinct. Assume first that they are, i.e. $\zeta_1\neq\zeta_2$. We
satisfy \thetag{7.1.9} by setting $\bar A^1_2=\zeta\bar A^2_2$ and
choosing $\bar A^2_1$ according to \thetag{8.2.3}. The equation for
$b'_1$ becomes \thetag{7.1.18}, and we can make $b'_1=0$ by choosing
either $\zeta=\zeta_1$ or $\zeta=\zeta_2$. We have to check that one
of these choices makes the transformation biholomorphic. By Lemma
8.1.3, this is equivalent to at least one of the roots $\zeta_1$,
$\zeta_2$ not being a root of
$$
p_2(\zeta)=a_1\zeta^2+2b_1\zeta+c_1.\tag8.5.2
$$ 
By comparing the coefficients of $p_1(\zeta)$ and $p_2(\zeta)$, we
find that $p_1(\zeta)$ and $p_2(\zeta)$ have the same set of roots if
and only if
$$
\left\{\aligned&b_1^2-a_1c_1=0\\&a_1c_2+a_2c_1=0.
\endaligned\right.\tag8.5.3
$$
However, the first equation of \thetag{8.5.3} is also equivalent to
$p_2(\zeta)$ having a 
double root. Since we assumed $\zeta_1\neq\zeta_2$, it follows that
both of these cannot be roots of $p_2(\zeta)$. 

Next, assume that $p_1(\zeta)$ has a double root. A straightforward
calculation shows that this happens if and only if
$$
(a_1c_2-a_2c_2)^2+4a_2c_2b_1^2=0,\tag8.5.5
$$
or equivalently,
$$
(a_1c_2+a_2c_1)^2+4a_2c_2(b_1^2-a_1c_1)=0.\tag8.5.6
$$
The double root is then
$$
\zeta_0=\frac{a_1c_2-a_2c_2}{2a_2b_1}.\tag8.5.7
$$
Now, \thetag{8.5.5} implies that $\zeta_0$ is also a root of $p_2(\zeta)$.
Thus, the transformation corresponding to the choice $\zeta=\zeta_0$ is
not biholomorphic. It follows that we cannot make $b'_1=0$. 
However, we can make $a'_2=0$ as follows. 

Assume first that $b^2_1-a_1c_1=0$, i.e. that $p_2(\zeta)$ also has a 
double root at $\zeta_0$. Then a straightforward calculation
shows that the transformation corresponding to $\bar A^2_2=0$, $\bar 
A^2_1=-\bar A^1_1a_1/b_1$, and $\bar A^1_2\neq0$ arbitrary, yields $a'_2 
=0$. 

Next, assume that $b^2_1-a_1c_1\neq0$. Equation \thetag{8.5.6} implies
then that $a_1c_2+a_2c_1\neq0$. The equation for $a'_2$ becomes
$$
\aligned
a'_2=&\,A^2_2(\bar A^1_1)^2\left(\bar\zeta\left(a_1-2b_1
\frac{s(\zeta,\bar\zeta)}{r(\zeta,\bar\zeta)}+c_1\left(
\frac{s(\zeta,\bar\zeta)}{r(\zeta,\bar\zeta)}\right)^2\right)
+a_2+c_2\left(
\frac{s(\zeta,\bar\zeta)}{r(\zeta,\bar\zeta)}\right)^2\right)\\=&\,
\frac{A^2_2(\bar A^1_1)^2}{r(\zeta,\bar\zeta)^2}
((a_1\bar\zeta+a_2)r(\zeta,\bar\zeta)^2-2b_1\bar\zeta
r(\zeta,\bar\zeta)s(\zeta,\bar\zeta)+(c_1\bar\zeta+c_2)
s(\zeta,\bar\zeta)^2)\\=&\,
\frac{A^2_2(\bar A^1_1)^2}{r(\zeta,\bar\zeta)^2}\cdot\\&\,
\left((c_1\bar\zeta+c_2)\left(s(\zeta,\bar\zeta)-\frac{b_1\bar\zeta
r(\zeta,\bar\zeta)}
{c_1\bar\zeta+c_2}\right)^2
+
\left(a_1\bar\zeta+a_2-
\frac{b_1^2\bar\zeta^2}{c_1\bar\zeta+c_2}\right)
r(\zeta,\bar\zeta)^2
\right)
\endaligned
\tag8.5.8
$$
Now, we solve for $b_1^2$ in \thetag{8.5.5} and obtain the following
$$
\aligned
a_1\bar\zeta+a_2-
\frac{b_1^2\bar\zeta^2}{c_1\bar\zeta+c_2}=&\,
a_1\bar\zeta+a_2+\frac{(a_1c_2-a_2c_1)^2}{4a_2c_2(c_1\bar\zeta+c_2)}
\bar\zeta^2
\\=&\,\frac
{4a_2c_2(a_1\bar\zeta+a_2)(c_1\bar\zeta+c_2)
+(a_1c_2-a_2c_1)^2\bar\zeta^2}{4a_2c_2(c_1\bar\zeta+c_2)}
\\=&\,
\frac
{(a_1c_2-a_2c_1)^2\bar\zeta^2+4a_2c_2(a_1c_2+a_2c_1)\bar\zeta+
4a_2^2c_2^2}{4a_2c_2(c_1\bar\zeta+c_2)}\\=&\,
\frac
{((a_1c_2+a_2c_1)\bar\zeta+2a_2c_2)^2}{4a_2c_2(c_1\bar\zeta+c_2)}.
\endaligned\tag8.5.9
$$
Moreover, using \thetag{7.1.16} and \thetag{8.5.6}, we obtain
$$
\aligned
s(\zeta,\bar\zeta)-\frac{b_1\bar\zeta}{c_1\bar\zeta+c_2}
r(\zeta,\bar\zeta)=&\,
a_1|\zeta|^2+b_1\bar\zeta+a_2\zeta-\frac{b_1\bar\zeta}{c_1\bar\zeta+c_2}
(b_1|\zeta|^2+c_1\bar\zeta+c_2)\\=&\,
\frac{(a_1c_1-b_1^2)\bar\zeta|\zeta|^2+(a_1c_2+a_2c_1)|\zeta|^2+a_2c_2\zeta}
{c_1\bar\zeta+c_2}\\=&\,
\frac{\zeta}{c_1\bar\zeta+c_2}\left(\frac{(a_1c_2+a_2c_1)^2}{4a_2c_2}\bar
\zeta^2+(a_1c_2+a_2c_2)\bar\zeta+a_2c_2\right)\\=&\,
\frac{\zeta}{4a_2c_2(c_1\bar\zeta+c_2)}
((a_1c_2+a_2c_1)\bar\zeta+2a_2c_2)^2.\endaligned\tag8.5.10
$$
If we combine \thetag{8.5.9} and \thetag{8.5.10}, then we get
$$
a'_2=\frac{A^2_2(\bar A^1_1)^2((a_1c_2+a_2c_1)\bar\zeta+2a_2c_2)^2}
{4a_2c_2(c_1\bar\zeta+c_2)
r(\zeta,\bar\zeta)^2} 
\left(\frac{\zeta^2((a_1c_2+a_2c_1)\bar\zeta+2a_2c_2)^2}{4a_2c_2}+
r(\zeta,\bar\zeta)^2\right).\tag8.5.11
$$
Consequently, if $\zeta_1$ is the solution of
$$
(a_1c_2+a_2c_1)\bar\zeta+2a_2c_2=0,\tag8.5.12
$$
and $\zeta_1$ is not a solution of $c_1\bar\zeta+c_2=0$ or 
\thetag{8.2.6}, then the transformation that corresponds
to the choice $\zeta=\zeta_1$ makes $a'_2=0$. The first condition
is easy to check,
$$
c_1\bar\zeta_1+c_2=c_2\frac{a_1c_2-a_1c_2}{a_1c_2+a_2c_1}\neq0\tag8.5.13
$$
Furthermore, 
$\zeta_1$ cannot be a solution of \thetag{8.2.6}, in view
of \thetag{8.5.10}. If it were, then $\zeta_1$ would also be a solution of
\thetag{8.2.17}. Lemma 8.1.1, in turn, would then imply that $\zeta_1$
is a root of $p_1(\zeta)$. This is not possible, since we assumed that
no root of $p_1(\zeta)$ solves \thetag{8.2.6}. 

In order to
conclude the proof of Theorem A (ii), we need to check that the
transformation corresponding to $\zeta=\zeta_1$ is biholomorphic.
A simple computation, using \thetag{8.5.10}, shows that 
$$
\bar A^1_1\bar A^2_2-\bar A^1_2\bar A^2_1=\bar A^1_1\bar A^2_2\frac{r(\zeta_1,\bar\zeta_1)^2}
{c_1\bar\zeta_1+c_2}\neq0.\tag8.5.14
$$
This completes the proof of Theorem A.

\heading 9. Proof of Theorem B\endheading
\subhead 9.1. The setup\endsubhead Following [CM], we shall reduce
the proof to a problem of describing the kernel and range of a certain
linear operator. We assume that $M$ is a hypersurface
in $\bC^3$ of the form (A.i.$k$), where $k$ is 1, 2, or 3. We write the defining equation of
$M$ at 0 as in \thetag{3.14}. We subject $M$ to a formal transformation
$$
z'=\tilde f(z,w)\quad,\quad w'=\tilde g(z,w),\tag9.1.1
$$
where $\tilde f=(\tilde f^1,\tilde f^2)$, 
that preserves the form of $M$ modulo terms of weighted
degree at least 4, i.e. the transformed hypersurface $M'$ is
given by a defining equation of the form 
$$
\im w' =|z'_1|^2+p_3(z',\bar z')+F'(z',\bar z',\re w'),\tag 9.1.2
$$
where $F'(z',\bar z', s')$ contains terms of weighted degree greater than
or equal to 4. We also require that the new
coordinates are regular for $M'$, i.e. that 
$$
\tilde g(z,w)=Q'\left(\tilde f(z,w),\overline{\tilde
f}(0,w),\overline{\tilde g}(0,w)\right),\tag9.1.3
$$
where $w'=Q'(z',\bar z',\bar w')$ is a complex defining equation of
$M'$ at 0 as in \S 5. Thus, $\tilde f$ and $\tilde g$ are subjected to the
restrictions imposed by Proposition 3.3. As mentioned in \S3, the most general
transformation of this kind can be factored uniquely as 
$$
(\tilde f(z,w),\tilde g(z,w))=(T\circ P)(z,w),\tag9.1.4
$$
where $P$ and $T$ are as described in \S3. 

To prove Theorem B, it suffices to prove that there is a unique
transformation 
$$
T(z,w)=(\hat f(z,w),\hat g(z,w))=(z+f(z,w),w+g(z,w))\tag 9.1.5
$$
to normal form (i.e. such that the transformed hypersurface $M'$ is
defined by \thetag{9.1.2} with $F'\in\Cal N^k$ for $k=1,2,3$
acccordingly) such that $f=(f^1,f^2)$, $f^1$ is $O(3)$, $f^2$
is $O(2)$, $g$ is $O(4)$, and such that the following additional
conditions are satisfied. If $(M,0)$ is of the form (A.i.1) or
(A.i.3), then the constant terms in the following formal series vanish
$$
\frac{\partial^2 f^2}{\partial z^\alpha}\,,\,  \re
\frac{\partial^2
f^1}{\partial z_1\partial w}\,,\,
\frac{\partial^3 f^1}{\partial z^\beta}.\tag9.1.6
$$
If $(M,0)$ is of the form (A.i.2), then the constant terms in the
following series vanish 
$$
\left\{\aligned
&\frac{\partial^2
f^j}{\partial z_k\partial w}\,,\,\frac{\partial^3 f^1}{\partial z^\beta}\,,\,\frac{\partial^3
f^1}{\partial z^\alpha \partial w}\,,\,\frac{\partial^3 f^j}{\partial
z_k\partial w^2}\\&
\frac{\partial^4 f^1}{\partial z^\beta\partial w}\,,\,\re
\frac{\partial^4
f^1}{\partial z_1\partial w^3}.\endaligned\right.\tag9.1.7
$$
We decompose $(f^1,f^2,g)$, $F$, and $F'$ into weighted homogeneous
parts as follows
$$
\aligned
f^1(z,w)=\sum_{\nu=3}^\infty f^1_\nu(z,w)\quad,\quad f^2(z,w) &=\sum_{\nu=2}^\infty
f^2_\nu(z,w)\quad,\quad g(z,w) =\sum_{\nu=4}^\infty g_\nu(z,w)\\
F(z,\bar z,s) =\sum_{\nu=4}^\infty F_\nu(z,\bar z,s)\quad&,\quad
F'(z,\bar z,s) =\sum_{\nu=4}^\infty F'_\nu(z,\bar z,s).\endaligned
$$
Recall here that $z$ and $\bar z$ are assigned the weight one, $w$ and
$s$ are assigned the weight two, and we say that e.g. $F_\nu(z,\bar
z,s)$ is weighted homogeneous of degree $\nu$ if for all $t>0$
$$
F_\nu(tz,t\bar z,t^2s)=t^\nu F_\nu(z,\bar z,s).
$$
The formal power series $F,F'\in\Cal F$ (recall, from \S3, the notation $\Cal F$
for the space of real power series with terms of weighted degree at
least 4) are related as follows 
$$
\aligned
\im \hat g(z,s+i\phi) &\equiv|\hat f^1(z,s+i\phi)|^2+p_3(\hat
f(z,s+i\phi),\overline {\hat f}(\bar
z,s-i\phi))\\& \,+F(\hat f(z,s+i\phi),\overline {\hat f}(\bar
z,s-i\phi),\re \hat g(z,s+i\phi)),\endaligned\tag9.1.8
$$
where
$$
\phi=\phi(z,\bar z,s)=|z_1|^2+p_3(z,\bar z)+F(z,\bar z,s).\tag9.1.9
$$
Identifying terms of weighted degree $\nu\geq 4$ we obtain
$$
F_\nu+\im
g_\nu\equiv z_1\overline{f^1_{\nu-1}}+\bar
z_1f^1_{\nu-1}+p_{3,z_2}f^2_{\nu-2}+p_{3,\bar
z_2}\overline{f^2_{\nu-2}}+F'_\nu+\ldots,\tag9.1.10
$$
where
$$
\aligned
F_\nu=F_\nu(z,\bar
z,s)\quad&,\quad F'_\nu=F'_\nu(z,\bar
z,s+i|z_1|^2)\\
p_{3,z_2}=p_{3,z_2}(z, &\bar z)=\frac{\partial p_3}{\partial z_3}(z,\bar
z)
\endaligned\tag9.1.11
$$
$$
\overline{f^1_{\nu-1}}=\overline{f^1_{\nu-1}}(\bar z,s-i|z_1|^2)\quad,\quad
f^1_{\nu-1}=f^1_{\nu-1}(z,s+i|z_1|^2)\quad,\quad\text{\rm etc},\tag9.1.12
$$
and where  the dots $\ldots$ signify terms that only involve $F_\mu$, $F'_\mu$, $g_\mu$,
$f^1_{\mu-1}$, and $f^2_{\mu-2}$ for $\mu<\nu$. Noting that 
$$
p_{3,\bar z_2}(z,\bar z)=\overline{p_{3,z_2}(z,\bar z)}\tag9.1.13
$$
we can write this as
$$
\re(ig_\nu+2\bar
z_1f^1_{\nu-1}+2p_{3,z_2}f^2_{\nu-2})=F_\nu-F'_\nu+\ldots.\tag9.1.14
$$
Let us define the linear operator 
$$
L(f^1,f^2,g)=\re(ig+2\bar
z_1f^1+2p_{3,z_2}f^2)|_{(z,s+i|z_1|^2)}\tag9.1.15
$$
from the space $\Cal G$ to the
space $\Cal F$, where $\Cal G$ denotes the space of formal power
series (in $(z,w)$) transformations 
$(f^1,f^2,g)$ 
such that $f^1$ is $O(3)$, $f^2$ is $O(2)$, and $g$ is $O(4)$. Observe
that $L$ maps $(f^1_{\nu-1},f^2_{\nu-2}, 
g_\nu)$ to a series that is weighted homogeneous of degree $\nu$. We
note, as in [CM], that if we could find subspaces
$$
\Cal G_0\subset\Cal G\quad,\quad\Cal N\subset\Cal F\tag9.1.16
$$
such that, for any $F\in\Cal F$, the equation
$$
L(f^1,f^2,g)=F\quad \mod \Cal N\tag9.1.17 
$$
has a unique solution $(f^1,f^2,g)\in\Cal G_0$ then, given any $F'\in\Cal F$, the
equation \thetag{9.1.14} would allow us to inductively determine 
the weighted homogeneous parts $F_\nu$ of a normal form $F\in
\Cal N$ and the weighted homogeneous parts
$(f^1_{\nu-1},f^2_{\nu-2},g_\nu)$ of the transformation
$(f^1,f^2,g)\in\Cal G_0$ to normal
form in a unique fashion. (This can
also be formulated as saying that $\Cal G_0$ and $\Cal N$ are
complementary subspaces of the kernel and range of $L$,
respectively). 

Let us therefore define $\Cal G_0^n\subset\Cal G$, for $n=1,2,3$, 
as those $(f^1,f^2,g)\in\Cal G$ for which the
constant terms, in the series \thetag{9.1.6} if $n=1,2$ or in
the series \thetag{9.1.7} if $n=2$, vanish. 
Combining all that is said above, we find that the proof of Theorem B
will be complete if we prove the following.

\proclaim{Lemma 9.1.18} Suppose the equation \thetag{9.1.2} is of the
form {\rm (A.i.$n$)} for some $n=1,2,3$.  Let $\Cal G^n_0\subset\Cal G$
be as described above and $\Cal N^n$ as defined in \S$3$. 
Then, for any $F\in\Cal F$, the equation
$$
L(f^1,f^2,g)=F\quad\mod\Cal N^n\tag9.1.19
$$
has a unique solution $(f^1,f^2,g)\in\Cal G^n_0$.\endproclaim
\demo{Proof} The proof is based on decomposing the equation
\thetag{9.1.19} according to type. Recall that a function
$G_{kl}(z,\bar z,s)$ has type $(k,l)$ if, for any $t_1,t_2>0$, the
following holds
$$
G_{kl}(t_1z,t_2\bar z,s)=t_1^kt_2^lG_{kl}(z,\bar z, s).
$$
We decompose $F\in \Cal F$ according to type
$$
F(z,\bar z,s)=\sum_{k,l}F_{kl}(z,\bar z,s).\tag9.1.20
$$
We also decompose $(f^1,f^2,g)\in\Cal G$ as follows ($j=1,2$)
$$
f^j(z,w)=\sum_{k}f^j_k(z,w)\quad,\quad
g(z,w)=\sum_{k}g_k(z,w),\tag9.1.21
$$
where $f^j_k(z,w)$, $g_k(z,w)$ are homogeneous of degree $k$ in $z$,
e.g.
$$
g_k(tz,w)=t^kg_k(z,w)\,,\quad t>0.\tag9.1.22
$$
The reader should observe that this redefines e.g. $g_2(z,w)$ which,
previously, denoted the weighted homogeneous part of degree $2$ in
$g(z,w)$. However, in what follows we shall not need the decomposition into weighted
homogeneous terms and, hence, the above notation should cause no
confusion; for the remainder of this paper, e.g. $g_2(z,w)$ means the part
of $g(z,w)$ which is homogeneous of degree 2 in $z$, etc. For brevity, we also use
the following notation (as in [CM])
$$
f'(z,w)=\frac{\partial f}{\partial
w}(z,w),\,\ldots,\,f^{(m)}(z,w)=\frac{\partial^m f}{\partial w^m}(z,w),\,\ldots.\tag9.1.23
$$
We will also use the fact
$$
f(z,s+i|z_1|^2)=f(z,s+iz_1\bar z_1)=\sum_{m}f^{(m)}(z,s)\frac{(iz_1\bar
z_1)^m}{m!}.\tag9.1.24
$$
We shall identify terms of type $(k,l)$ in \thetag{9.1.19}. Since the
equation is real, it suffices to consider types where $k\geq l$. Also,
note that for $(k,l)$ such that $\Cal N^n_{kl}=\Cal F_{kl}$ the
equation \thetag{9.1.19} is trivially satisfied. (We use the
obvious notation; e.g. $\Cal F_{kl}$ are the elements in $\Cal F$ that
are of type $(k,l)$, etc.) We will need to decompose $p_{3,z_2}(z,\bar z)$
according to type. We use the notation $p_{11}$ for the part
of $p_{3,z_2}$ which is of type $(1,1)$ and $p_{02}$ for the part
which is if type $(0,2)$. Thus,
$$
p_{3,z_2}(z,\bar z)=p_{11}(z,\bar z)+p_{02}(z,\bar z)\quad,\quad
p_{3,\bar z_2}(z,\bar z)=\overline{p_{11}(z,\bar
z)}+\overline{p_{02}(z,\bar z)}.\tag9.1.25
$$
Here, the cases (A.i.1) and (A.i.3) differ from (A.i.2) in the respect
that $p_{11}(z,\bar z)\equiv 0$ in the latter case whereas
$p_{11}(z,\bar z)\not\equiv 0$ in the first two. This fact makes the
normal form for (A.i.2) different from the normal forms for (A.i.1)
and (A.i.3). The conclusion of the proof will therefore be carried out in two
versions, one for (A.i.1), (A.i.3) and one for (A.i.2).

\enddemo

\subhead 9.2. Conclusion of the proof of Lemma 9.1.18; the cases
(A.i.1) and (A.i.3)\endsubhead In what follows, we use the notation
$$
F_{kl}=F_{kl}(z,\bar z,s)\quad,\quad g_k=g_k(z,s)\quad,\quad
\overline{g_k}=\overline{g_k}(\bar z,s) 
\quad,\quad\text{\rm
etc}.\tag9.2.1
$$
Collecting terms of equal type in
\thetag{9.1.19}, we obtain the following decoupled systems of differential
equations, for $k\geq 3$,
$$
\left\{
\aligned
\frac{i}{2}g_k &=F_{k0}\\
\bar z_1f^1_{k+1}+p_{11}f^2_k-\frac{z_1\bar z_1}{2}g'_k
&=F_{k+1,1}\,\mod\Cal N^n_{k+1,1},
\endaligned
\right.
\tag9.2.2
$$
where $n=1,3$ according to whether the type is (A.i.1) or (A.i.3), and
$$
\align
&\left\{
\aligned
\pb02\overline{f^2_0}+\frac{i}{2}g_2 &=F_{20}\\
\bar z_1f^1_3+p_{11}f^2_2-iz_1\bar z_1\pb02(\overline{f^2_0})'-
\frac{z_1\bar z_1}{2}g'_2 &=F_{31}\,\mod \Cal N_{31}^n\\
iz_1\bar z_1^2(f^1_3)'+iz_1\bar z_1p_{11}(f^2_2)'-\frac{z_1^2\bar
z_1^2}{2}\pb02 (\overline{f^2_0})''&\\+p_{02}f^2_4-\frac{iz_1^2\bar
z_1^2}{4}g''_2 &=F_{42}
\,\mod\Cal N^n_{42},\endaligned
\right.
\tag9.2.3
\\
&\left\{
\aligned
z_1\overline{f^1_0}+\frac{i}{2}g_1 &=F_{10}\\
-iz_1^2\bar z_1(\overline{f^1_0})'+\bar
z_1f^1_2+p_{11}f^2_1+\pb02\overline{f^2_1}-\frac{z_1\bar z_1}{2}g'_1
&=F_{21}\,\mod\Cal N^n_{21}\\
iz_1\bar z_1^2(f^1_2)'-\frac{z_1^3\bar
z_1^2}{2}(\overline{f^1_0})'' +iz_1\bar
z_1p_{11}(f^2_1)'&\\+p_{02}f^2_3-iz_1\bar
z_1\pb02(\overline{f^2_1})'-\frac{iz_1^2\bar z_1^2}{4}g_1''
&=F_{32}\,\mod\Cal N^n_{32}\endaligned
\right.
\tag9.2.4
\\
&\left\{
\aligned
-\im g_0 &=F_{00}\\
2\re(\bar z_1f^1_1)+2\re(p_{11}f^2_0)-z_1\bar z_1\re
g'_0 &=F_{11}\,\mod\Cal N^m_{11}\\
2z_1\bar z_1\im(p_{11}(f^2_0)')-2\re(p_{02}f^2_2)&\\-2z_1\bar z_1\im(\bar
z_1(f^1_1)')+\frac{z_1^2\bar z_1^2}{2}\im g''_0 &=F_{22}\,\mod\Cal
N^n_{22}\\
-z_1^2\bar z_1^2\re(\bar z_1(f^1_1)'')-z_1^2\bar
z_1^2\re(p_{11}(f^2_0)'')&\\-2z_1\bar
z_1\im(p_{02}(f^2_2)')+\frac{z_1^3\bar z_1^3}{6}\re
g'''_0 &=F_{33}\,\mod\Cal N^n_{33}.\endaligned
\right.
\tag9.2.5
\endalign
$$
Solving \thetag{9.2.2} is easy. Substitute $g_k$ in the second
equation
$$
\bar z_1f^1_{k+1}+p_{11}f^2_k=F_{k+1,1}+\ldots\,\mod \Cal
N^n_{k+1,1},\tag9.2.6
$$
where $\ldots$ signifies known terms whose precise form is not important
since we are solving mod $\Cal N^n_{k+1,1}$ (in this case though,
$\ldots$ is easy to compute). Now, when $n=1$ we have
$$
p_{11}(z,\bar z)=2z_2\bar z_2\tag9.2.7
$$
and when $n=3$
$$
p_{11}(z,\bar z)=z_1\bar z_2.\tag9.2.8
$$
Thus, in both cases we can find $f^1_{k+1}$ and $f^2_k$ uniquely by
decomposing the right side as $\bar z_1 G_{k+1,0}+\bar z_2
H_{k+1,0}$. However, since $p_{11}$ also contains a factor $z_j$, where
$j$ is 1 or 2 according to whether $n$ is 1 or 3, we can only solve
\thetag{9.2.6} if $\Cal N^n_{k+1,1}$ is a complemetary subspace of those $F_{k+1,1}\in\Cal
F_{k+1,1}$ that are of the form
$$
F_{k+1,1}=\bar z_1 G_{k+1,0}+\bar z_2z_jH_{k0}.\tag9.2.9
$$
Clearly, $\Cal N^n_{k+1,1}$ described in \thetag{3.25} is such a
subspace. Hence, \thetag{9.2.6} can be solved uniquely for
$f^1_{k+1}$, $f^2_k$, and $g_k$, for all $k\geq 3$.

We proceed to solve \thetag{9.2.3}. Solve for $g_2$ in the first
equation and substitute the result in the last two. Noting in the
third equation that the 
terms involving $\overline{f^2_0}$ cancel and that the
function $f^2_4$ is known from the previous step,
this third equation equation becomes 
$$
iz_1\bar z_1^2(f^1_3)'+iz_1\bar
z_1p_{11}(f^2_2)'=F_{42}+\ldots\mod\Cal N^n_{42}.\tag9.2.10
$$
As above, the fact that $p_{11}$ is divisible by $\bar z_2$ implies
that \thetag{9.2.10} determines uniquely $(f^1_3)'$ and
$(f^2_2)'$. Thus, $f^1_3$ and $f^2_2$ are uniquely determined given a
choice of 
$$
f^1_3(z,0)\quad,\quad f^2_2(z,0).\tag9.2.11
$$
Since $f^1_3(z,0)$ and $f^2_2(z,0)$ denote the parts of $f^1(z,0)$ and
$f^2(z,0)$ that are homogeneous of degrees 3 and 2, respectively, the
choice in \thetag{9.2.11} amounts to a choice of constant terms in the
series
$$
\frac{\partial^3 f^1}{\partial
z^\beta}\quad,\quad\frac{\partial^2f^2}{\partial z^\alpha},\tag9.2.12
$$
where $\alpha$ ranges over all multi-indices with $|\alpha|=2$ and
$\beta$ over those with $|\beta|=3$. Thus, $f^1_3$ and $f^2_2$ are
uniquely determined by the condition that $(f^1,f^2,g)\in\Cal G^n_0$. 
In order for \thetag{9.2.10} to
be satisfied with this choice of $f^1_3$ and $f^2_2$, $\Cal N^n_{42}$
has to be complementary to those $F_{42}\in\Cal F_{42}$ that are of
the form described by the left hand side of \thetag{9.2.10}. It is
straightforward to verify that $\Cal N^n_{42}$ described by
\thetag{3.24} is such. Substitute $f^1_3$ and $f^2_2$ so obtained in
the second equation. We obtain
$$
iz_1\bar z_1\pb02(\overline{f^2_0})'=F_{31}+\ldots\mod\Cal N^n_{31}.\tag
9.2.13
$$
This determines $f^2_0$ uniquely since the constant term in the
series $f^2$ must be 0. Also, it is easy to verify, using the fact
that
$$
\aligned
\pb02(\bar z,z) =z_2^2+\gamma z_1^2\quad &,\quad\text{\rm if $n=1$}\\
\pb02(\bar z,z) =z_1z_2\quad &,\quad\text{\rm if $n=3$},\endaligned
\tag9.2.14
$$
that then
\thetag{9.2.13} is satisfied if $\Cal N^n_{31}$ is defined by
\thetag{3.22}. 

Next, we turn to the system \thetag{9.2.4}. Solve for $g_1$ in the
first equation and substitute into the last two. Then differentiate
the second equation, multiply it by $iz_1\bar z_1$, and subtract from
the third equation. We obtain, noting that the series $f^2_3$ was
determined from the system \thetag{9.2.2}, 
$$
-2iz_1\bar z_1\pb02(\overline{f^2_1})'-2z_1^3\bar
z_1^2(\overline{f^1_0})''=F_{32}+\ldots\mod\Cal N^n_{32}.\tag9.2.15
$$
Using \thetag{9.2.14}, we deduce that \thetag{9.2.15} uniquely
determines $f^2_1$, since the constant term in 
$\partial f^2/\partial z_j$ must be 0 for $j=1,2$. Then, the same equation
determines $f^1_0$ uniquely; the constant terms in $f^1$ and 
$\partial f^1/\partial w$ must be 0. With
$f^2_1$ and $f^1_0$ defined this way, \thetag{9.2.15} is satisfied if
$\Cal N^n_{32}$ is defined by \thetag{3.23}. Now, substitute the
series so obtained back into the second equation and find
$$
\bar z_1f^1_2=F_{21}+\ldots\mod\Cal N^n_{21}.\tag9.2.16
$$
This determines $f^1_2$ uniquely, and \thetag{9.2.16} is satified with
$\Cal N^n_{21}$ defined by \thetag{3.21}.

The last step is to solve the system of equations
\thetag{9.2.5}. Substituting the first equation in the third and
observing that $f^2_0$ and $f^2_2$ have already been determined, we
obtain the following
$$
-2z_1\bar z_1\im(\bar z_1(f^1_1)')=F_{22}+\ldots\mod\Cal
N^n_{22}.\tag9.2.17
$$
Now, the series $f^1_1$ can be written
$$
f^1_1(z,s)=a_1(s)z_1+a_2(s)z_2\tag 9.2.18
$$
and, hence,
$$
\im(\bar z_1(f^1_1)')=\im a'_1z_1\bar z_1+\frac{1}{2i}a'_2z_2\bar
z_1-\frac{1}{2i}\overline{a'_2}z_1\bar z_2.\tag 9.2.19
$$
It follows that $\im a_1$ and $a_2$ can be uniquely determined from
\thetag{9.2.17}, since
the constant terms in $\im a_1$ and $a_2$ are 0. Also, the equation
\thetag{9.2.17} is satisfied if 
$\Cal N^n_{22}$ is defined as in \thetag{3.21}. Next, differentiate
the second equation of \thetag{9.2.5} twice, multiply by $z_1^2\bar z_1^2/6$, and add to
the fourth. We obtain
$$
-\frac{2z_1^2\bar z_1^2}{3}\re(\bar
z_1(f^1_1)'')=F_{33}+\ldots\mod\Cal N^n_{33}.\tag9.2.20
$$
Using \thetag{9.1.18}, we can write
$$
\re(\bar z_1(f^1_1)'')=\re a''_1z_1\bar z_1+\frac{1}{2}a''_2z_2\bar
z_1+\frac{1}{2}\overline{a''_2}z_1\bar z_2.\tag9.2.21
$$
Hence, \thetag{9.2.20} determines $\re a_1$ uniquely, with a choice of
constant term in $\re a'_1$, and \thetag{9.2.20} is
then satisfied if $\Cal N^n_{33}$ is as defined as in \thetag{3.21}.
We conclude that $f^1_1$ is uniquely determined with a choice of
constant term in 
$$
\re\frac{\partial^2f^1}{\partial z_1\partial w}.\tag
9.2.22
$$
Hence, $f^1_1$ is uniquely determined by the condition
$(f^1,f^2,g)\in\Cal G^n_0$. 
Substituting back into the second equation, we obtain
$$
-z_1\bar z_1\re g'_0=F_{11}+\ldots\mod\Cal N^n_{11}.\tag9.2.23
$$
We deduce that $\re g_0$ is uniquely determined, since the constant
term in $\re g_0$ must be 0, and \thetag{9.2.23} is satisfied if $\Cal
N^n_{11}$ is given by \thetag{3.21}. This completes the proof of Lemma
9.1.18 in the cases (A.i.$n$) with $n$ either 1 or 3.

\subhead 9.3. Conclusion of the proof of Lemma 9.1.18; the case
(A.i.2)\endsubhead The main difference between this situation and the
previous one is, as already mentioned above, that $p_{11}=0$. This
implies e.g. that the system \thetag{9.2.2} does not involve $f^2_k$, the
system \thetag{9.2.3} does not involve $f^2_2$. Consequently, we have
to add more equations. We shall add the following equations that arise
from \thetag{9.1.19}; in what follows,
since we are now dealing with just one case, we shall make the
replacement 
$$
p_{02}(z,\bar z)=\bar z_1^2.\tag9.3.1
$$
To the system \thetag{9.2.2}, we add the equation
$$
iz_1\bar z_1^2(f^1_{k+1})' +\bar z_1^2f^2_{k+2}-\frac{iz_1^2\bar
z_1^2}{4}g_k''=F_{k+2,2}\mod \Cal N^3_{k+2,1},\tag9.3.2
$$
for $k\geq 3$. To the system \thetag{9.2.3} we add
$$
-\frac{z_1^2\bar z_1^3}{2}(f^1_3)''+iz_1\bar
z_1^3(f^2_4)'+\frac{iz_1^5\bar
z_1^3}{6}(\overline{f^2_0})'''+\frac{z_1^3\bar
z_1^3}{12}g'''_2=F_{53}\mod\Cal N^3_{53}.\tag9.3.3
$$
To the systems \thetag{9.2.4} and \thetag{9.2.5}, we add respectively
$$
\align
&\left\{
\aligned
-\frac{z_1^2\bar z_1^3}{2}(f^1_2)''+\frac{iz_1^4\bar z_1^3}{6}
(\overline{f^1_0})'''+iz_1\bar z_1^3(f^2_3)'&\\-\frac{z_1^4\bar
z_1^2}{2}(\overline{f^2_1})''+\frac{z_1^3\bar
z_1^3}{12}g'''_1 &=F_{43}\mod\Cal N^3_{43}\\
-\frac{iz_1^3\bar z_1^4}{6}(f^1_2)'''+\frac{z_1^5\bar
z_1^4}{24}(\overline{f^1_0})^{(4)}-\frac{z_1^2\bar
z_1^4}{2}(f^2_3)''&\\+\frac{iz_1^5\bar
z_1^3}{6}(\overline{f^2_1})'''+\frac{iz_1^4\bar
z_1^4}{48}g_1^{(4)} &=F_{54}\mod\Cal
N^3_{54}.\endaligned\right.\tag9.3.4
\\
&\left\{
\aligned
\frac{z^3_1\bar z_1^3}{3}\im(\bar z_1(f^1_1)''')-z_1^2\bar
z_1^2\re(\bar z_1^2(f^2_2)'')-\frac{z_1^4\bar z_1^4}{24}\im g^{(4)}_0
&=F_{44}\mod\Cal N^3_{44}\\
\frac{z_1^4\bar z_1^4}{12}\re(\bar z_1(f^1_1)^{(4)})+\frac{z_1^3\bar
z_1^3}{3}\im (\bar z_1^2(f^2_2)''')-\frac{z_1^5\bar
z_1^5}{120}\re g_0^{(5)} &=F_{55}\mod\Cal N^3_{55}.
\endaligned\right.\tag9.3.5
\endalign
$$
The systems of differential equations so obtained (by adding the above
to the systems in \S9.2) can be solved if we define the space $\Cal
N^3$ as in \S3, and the solution $(f^1,f^2,g)$ is unique if we require
it to be in $\Cal G^2_0$ as defined in the paragraph preceding Lemma
9.1.18. The proof of 
this proceeds along the same lines as the proof in the
previous section although, due to the fact that these systems have
more equations, the calculations become a bit more tedious. The
modifications, however, are straightforward and the details are left
to the reader. As mentioned above, this also 
completes the proof of Theorem B.\qed


\Refs\widestnumber\key{BER2}

\ref\key BER1 \manyby M. S. Baouendi, P. Ebenfelt and
L. P. Rothschild\paper Algebraicity of holomorphic mappings between
real algebraic sets in $\bC^N$
\jour Acta Math.
\vol 177\yr 1996\pages 225--273
\endref

\ref\key BER2 \bysame\paper Infinitesimal CR
automorphisms of real analytic manifolds in complex
space
\jour Comm. Anal. Geom.
\finalinfo(to appear)
\endref

\ref\key BER3 \bysame\paper Parametrization of local biholomorphisms
of real analytic hypersurfaces
\jour Asian J. Math.
\finalinfo(to appear)
\endref

\ref\key BER4 \bysame\book Real manifolds in complex space
\finalinfo(in preparation; preliminary title)
\endref


\ref\key BJT\by M. S. Baouendi, H. Jacobowitz and F.
Treves\paper On the analyticity of CR mappings\jour
Ann. of Math.\vol 122\yr 1985\pages 365--400
\endref


%
%
%
%
%

%
%

\ref\key BS\by D. Burns, Jr. and S. Shnider\paper Real hypersurfaces
in complex manifolds\inbook Proceedings of
Symposia in Pure Mathematics XXX, Part 2, Several Complex Variables\publaddr Amer
Math. Soc., Providence, RI\yr 1977\pages 141--168\endref


\ref\key C1\manyby E. Cartan\paper Sur la g\'eom\'etrie
pseudo-conforme des hypersurfaces de deux variables complexes, I\jour
Ann. Math. Pura Appl.\vol 11\yr 1932\pages 17--90\finalinfo (or Oeuvres
II, 1231--1304)\endref

\ref\key C2\bysame \paper Sur la g\'eom\'etrie
pseudo-conforme des hypersurfaces de deux variables complexes, II\jour
Ann. Scoula Norm. Sup. Pisa\vol 1\yr 1932\pages 333--354\finalinfo(or
Oeuvres III, 1217--1238)\endref

\ref\key CM
\by S.-S. Chern and J.K. Moser
\paper Real hypersurfaces in complex manifolds
\jour Acta Math.
\vol 133 \yr 1974 \pages 219-271
\endref

%
%

\ref\key F\by M. Freeman\paper Real submanifolds with degenerate Levi
form\inbook Proceedings of Symposia in Pure Mathematics XXX, part I, Several
Complex Variables\publ Amer. Math. Soc.\publaddr Providence, RI\yr 1977\endref

%
%
%
%
%
%

\ref\key P
\by  H. Poincar\'e
\paper Les fonctions analytiques de deux variables et la repr\'esentation
conforme
\jour Rend. Circ. Mat. Palermo, II.  Ser. 
 \vol 23 \yr 1907 \pages 185-220
\endref

%
%
%

\ref\key S
\manyby  N. Stanton
\paper Infinitesimal CR automorphisms of rigid hypersurfaces
\jour Amer. J. Math. 
\vol 117 \yr 1995
\pages 141-167
\endref

%
%
\ref\key T1
\manyby  N. Tanaka
\paper On the pseudo-conformal geometry of hypersurfaces of the space of
$n$ complex variables
\jour J. Math. Soc. Japan 
\vol 14 \yr 1962 \pages 397-429
\endref

\ref\key T2
\bysame
\paper On generalized graded Lie algebras and geometric
structures. I\jour J. Math. Soc. Japan \vol 19\yr 1967\pages
215--254\finalinfo (erratum {\bf 36}, p. 1568)\endref


%
%
%
%
\ref\key W\by S. M. Webster\paper The holomorphic contact geometry of
a real hypersurface\inbook Modern Methods in Complex Analysis\eds
T. Bloom et al\publ Annals of Mathematics Studies 137, Princeton
University Press\publaddr Princeton, N.J.\yr 1995\pages 327--342\endref
%

%
%
%

\endRefs
\enddocument
\end